\documentclass[12pt]{article}

\usepackage{latexsym}
\usepackage{amsfonts}
\usepackage{amssymb}
\usepackage{amsmath}
\usepackage{mathrsfs}
\usepackage{bbm}
\usepackage{dsfont}
\usepackage{hyperref}

\newcommand{\be}{\begin{equation}}
\newcommand{\ee}{\end{equation}}
\newcommand{\bea}{\begin{eqnarray}}
\newcommand{\eea}{\end{eqnarray}}
\newcommand{\bean}{\begin{eqnarray*}}
\newcommand{\eean}{\end{eqnarray*}}
\newcommand{\brray}{\begin{array}}
\newcommand{\erray}{\end{array}}

\newtheorem{dfn}{Definition}[section]
\newtheorem{thm}[dfn]{Theorem}
\newtheorem{lmma}[dfn]{Lemma}
\newtheorem{ppsn}[dfn]{Proposition}
\newtheorem{crlre}[dfn]{Corollary}
\newtheorem{xmpl}[dfn]{Example}
\newtheorem{rmrk}[dfn]{Remark}

\newcommand{\bdfn}{\begin{dfn}\rm}
\newcommand{\bthm}{\begin{thm}}
\newcommand{\blmma}{\begin{lmma}}
\newcommand{\bppsn}{\begin{ppsn}}
\newcommand{\bcrlre}{\begin{crlre}}
\newcommand{\bxmpl}{\begin{xmpl}}
\newcommand{\brmrk}{\begin{rmrk}\rm}

\newcommand{\edfn}{\end{dfn}}
\newcommand{\ethm}{\end{thm}}
\newcommand{\elmma}{\end{lmma}}
\newcommand{\eppsn}{\end{ppsn}}
\newcommand{\ecrlre}{\end{crlre}}
\newcommand{\exmpl}{\end{xmpl}}
\newcommand{\ermrk}{\end{rmrk}}

\newcommand{\clh}{\mathcal{H}}

\newcommand{\clk}{\mathcal{K}}

\author{ S. Sundar}
\title {Arveson's characterisation of CCR flows : the multiparameter case}

\begin{document}
\maketitle
\begin{abstract}
In this paper, we revisit Arveson's characterisation of CCR flows in terms of decomposibility of the product system  in the multiparameter context. We show that a multiparameter $E_0$-semigroup is a CCR flow if and only if it is decomposable and admits a unit. In contrast to the one parameter situtation, we exhibit  uncountably many examples of  decomposable $E_0$-semigroups which do not admit any unit. As applications, we show that for a pure isometric representation $V$, the associated CCR flow $\alpha^{V}$ remembers the unitary equivalence class of $V$. We also compute the positive contractive local cocycles and the projective local cocycles of a CCR flow. A necessary and a sufficient condition for a CCR flow to be prime is obtained. 
 \end{abstract}
\noindent {\bf AMS Classification No. :} {Primary 46L55; Secondary 46L99.}  


\section{Introduction}
The theory of $E_0$-semigroups initiated by R.T.Powers and further developed extensively by  William Arveson in his seminal papers \cite{Arv_Fock}, \cite{Arv_Fock3}, \cite{Arv_Fock2} and \cite{Arv_Fock4} has been an active area of research for the past three decades. For Powers' influential work on the subject, we refer the reader to \cite{Powers_TypeIII}, \cite{Powers_Index} and \cite{Powers_CPflow}. For a more comprehensive list of references on the subject and a thorough treatment,  the reader is referred to the monograph \cite{Arveson}. Let $\clh$ be an infinite dimensional separable Hilbert space and $B(\clh)$ be the $*$-algebra of bounded operators on $\clh$. An $E_0$-semigroup on $B(\clh)$ is a $1$-parameter semigroup $\{\alpha_{t}\}_{t \geq 0}$ of unital normal $*$-endomorphisms of $B(\clh)$ which satisfies a natural continuity hypothesis. 

However, mathematically speaking, there is no reason to restrict our attention to  semigroups of endomorphisms indexed by the positive real line. Recently the author in collaboration with others has studied  semigroups of endomorphisms on $B(\clh)$ where the indexing set is a closed convex cone in a Euclidean space. The relationship between such semigroups and the associated product systems were explored in \cite{Murugan_Sundar} and in \cite{Murugan_Sundar_continuous}. Several authors have tried to understand noncommutative dynamics over several variables. The notable papers in this direction are \cite{Shalit_2008}, \cite{Shalit}, \cite{Shalit_2011}, \cite{Shalit_Solel} , \cite{Solel},  and \cite{Hirshberg_Daniel}.

Let us recall the definition of an $E_0$-semigroup over a closed convex cone. Fix a closed convex cone  $P$ in $ \mathbb{R}^{d}$ which we assume is  spanning and pointed, i.e. $P-P=\mathbb{R}^{d}$ and $P \cap -P=\{0\}$. Let $\clh$ be an infinite dimensional separable Hilbert space. By an $E_0$-semigroup over $P$ on $B(\clh)$, we mean a family $\alpha:=\{\alpha_{x}\}_{x \in P}$ of unital normal $*$-endomorphisms of $B(\clh)$ such that the following conditions hold: 
\begin{enumerate}
\item[(1)] for $x,y \in P$, $\alpha_{x}\circ \alpha_{y}=\alpha_{x+y}$, and
\item[(2)] given $A \in B(\clh)$ and $\xi,\eta \in \clh$, the map $P \ni x \to \langle \alpha_{x}(A)\xi|\eta \rangle \in \mathbb{C}$ is continuous.
\end{enumerate}
Since $P$ will be a fixed but an arbitrary cone, we will drop the phrase ``over $P$" and call our objects simply $E_0$-semigroups.

The simplest examples of $E_0$-semigroups arise out of the CCR construction. Let $V:=\{V_{x}\}_{x \in P}$ be a strongly continuous semigroup of isometries on a separable Hilbert space $\clh$. We also call such semigroups of isometries as  isometric representations of $P$. Denote the symmetric Fock space of $\clh$ by $\Gamma(\clh)$. Then there exists a unique $E_0$-semigroup $\alpha^{V}:=\{\alpha_{x}\}_{x \in P}$ on $B(\Gamma(\clh))$ such that for each $x \in P$ and $\xi \in \clh$,
\[
\alpha_{x}(W(\xi))=W(V_{x}\xi).\]
Here $\{W(\xi):\xi \in \clh\}$ denotes the usual Weyl operators. Recall that the action of the Weyl operators on the exponential vectors $\{e(\eta): \eta \in \clh\}$ is given by the formula
\[
W(\xi)e(\eta)=e^{-\frac{||\xi||^{2}}{2}-\langle \eta|\xi \rangle}e(\xi+\eta).\]
We call an isometric representation $V=\{V_{x}\}_{x \in P}$ pure if $\displaystyle \bigcap_{x \in P}V_{x}\clh=\{0\}$.

Consider the case when $P=[0,\infty)$. Arveson's fundamental work states that the map $V \to \alpha^{V}$ sets up a bijection between the set of pure isometric representations (up to unitary equivalence)
and the set of decomposable $E_0$-semigroups (up to cocycle conjugacy). 
\begin{enumerate}
\item[(1)] The proof of the injectivity part, undertaken in \cite{Arv_Fock}, relies heavily on the fact that there are enough units and on an index computation. 
\item[(2)] The proof of the surjectivity part, carried out in \cite{Arv_Path}, is far deeper. Arveson proves his result through a path space construction after surmounting difficult cohomological problems.
\end{enumerate}

In this paper, we address the injectivity  and the surjectivity question in the higher dimensional case. The injectivity question was analysed for a subclass of isometric representations that arise out of shifts in \cite{Anbu_Sundar}. The analysis carried out in \cite{Anbu_Sundar} makes heavy use of groupoid techniques and the analysis was possible precisely because the isometric representations considered in \cite{Anbu_Sundar} were shifts. 

One of the first difficulties in imitating Arveson's proof of  injectivity is that, in the multiparameter case,  the examples that we know so far  admits only one unit (up to a character). Consequently index computation does not reveal anything significant. Secondly, unlike in the one parameter case, there is no known Wold decomposition result in the multiparameter case, i.e. there is no good coordinatization of an isometric representation. Thus one is forced to look for a  coordinate free proof. 

Fortunately, Arvesons's ideas in \cite{Arv_Path} can indeed be turned around to yield  such a proof. In \cite{Arv_Path}, a strategy to construct an isometric representation from a decomposable product system is given. We imitate this construction in the higher dimensional setting. The key construction, whose construction forms the heart of this paper,  is to construct Arveson's e-logarithm in this setting which we achieve by appealing to the geometry of the cone.  With the $e$-logarithm in hand,  Arveson's arguments in the $1$-parameter setting carries over and characterise CCR flows as precisely those $E_0$-semigroups which are  decomposable and admit a unit. Despite the fact that we imitate Arveson, the author believes that  the construction of Arveson's $e$-logarithm in the higher dimensional setting and the applications given are non-trivial and are worth recording. 

We obtain the injectivity of the map $V \to \alpha^{V}$ as a byproduct of our construction. We show that multiparameter CCR flows are decomposable in a suitable sense and when we apply Arveson's construction to the CCR flow $\alpha^{V}$, we get back the isometric representation $V$ thereby proving the injectivity of the map $V \to \alpha^{V}$. This involves determining the decomposable vectors of the usual $1$-parameter CCR flows which we determine in Section 2. Prof. Liebscher has indicated to the author via e-mail that he believes that this is probably well known to experts. We include the details for completeness. 
The organisation of this paper is as follows. 

After determining the decomposable vectors of one parameter CCR flows in Section 2, we take up the construction of the $e$-logarithm in Section 3. 
The characterisation of CCR flows as decomposable $E_0$-semigroups admitting a unit is explained in Section 4. Uncountably many examples of   decomposable $E_0$-semigroups which are unitless are constructed.  In Section 5, we prove that the  map $V \to \alpha^{V}$ is injective. In Sections 6 and 7, we  give a few applications. In particular, we explain how the results obtained here provide a more conceptual explanation for the results obtained in \cite{Anbu} and \cite{Anbu_Sundar}. We also work out the  positive contractive local cocycles  and the projective local cocycles of CCR flows. The positive contractive local cocycles for $1$-parameter CCR flows were computed by Bhat in \cite{Bhat_cocycles}. In Section 7, we derive a necessary and a sufficient condition for a  CCR flow to be prime which means that it cannot be written as a tensor product of two $E_0$-semigroups. We show that $\alpha^{V}$ is prime if and only if the isometric representation $V$ is irreducible.

We assume that the reader is familiar with the terminology of Chapters 5 and 6 of \cite{Arveson} which we use extensively without recalling them.  We assume tacitly that all the Hilbert spaces involved are over $\mathbb{C}$ and are separable. Moreover our convention is that the inner product is linear in the first variable and conjugate linear in the second variable.

\section{Decomposable vectors of CCR flows}
In this section, we work out the set of decomposable vectors of a $1$-parameter CCR flow. 
Fix a one parameter $E_0$-semigroup  $\alpha:=\{\alpha_{t}\}_{t \geq 0}$ on $B(\clh)$.
Denote the product system associated to $\alpha$ by $E$ and the fibre of $E$ at a point $t$ will be denoted by $E(t)$. The set of decomposable vectors at $t$ is denoted by $D(t)$. Let \[\displaystyle \Delta:=\coprod_{t >0}\Delta(t)\] be the path space associated to $E$. Let us recall the definition of $\Delta(t)$. Fix $t>0$. For $x,y \in D(t)$, we say $x \sim y$ if and only if there exists $\lambda \in \mathbb{C}\backslash\{0\}$ such that $x=\lambda y$. Then $\sim$ is an equivalence relation on $D(t)$ and $\Delta(t)$ is the set of equivalence classes. For $x \in D(t)$, we denote the equivalence class of $x$ by $\dot{x}$. Set \[
\Delta^{(2)}:=\{(t,\dot{x},\dot{y}):t>0,x,y \in D(t)\}.\]

One of the key construction in \cite{Arv_Path} is the construction of the $e$-logarithm which is a positive definite function on the set of decomposable vectors.  Let $\displaystyle e:=(e_{t})_{t>0} \in \coprod_{t>0}D(t)$ be a left coherent section of decomposable vectors such that $||e_{t}||=1$. This means that for $s<t$, there exists a necessarily unique element $e(s,t) \in D(t-s)$ such that $e_{t}=e_{s}e(s,t)$.

Arveson's theorem is that there exists a unique continuous function $L^{e}:\Delta^{(2)} \to \mathbb{C}$ which vanishes at the origin and satisfies the following equation:
\[
e^{L^{e}(t,\dot{x},\dot{y})}=\frac{\langle x|y \rangle}{\langle x|e_{t}\rangle \langle e_{t}|y \rangle}
\]
for $t>0$, $x, y \in D(t)$. Moreover for every $t>0$, the map \[D(t) \times D(t) \ni (x,y) \to L^{e}(t,\dot{x},\dot{y}) \in \mathbb{C}\] is positive definite.

\begin{rmrk}
\begin{enumerate}
\item[(1)] For the precise definition of continuity of functions defined on $\Delta^{(2)}$ and $\Delta$, we refer the reader to the second paragraph of Section 6.4, Chapter 6 in \cite{Arveson}.
\item[(2)] For $x,y \in D(t)$ and $t>0$, we often abuse notation and write $L^{e}(t,x,y)$ instead of  $L^{e}(t,\dot{x},\dot{y})$. Such abuse of notation,  
while dealing with functions defined on $\Delta$ and $\Delta^{(2)}$, occur throughout this paper.  
\end{enumerate}
\end{rmrk}

Let $V:=(V_{t})_{t \geq 0}$ be a $1$-parameter isometric representation on a separable Hilbert space $\clh$. Denote the corresponding CCR flow on $B(\Gamma(\clh))$ by $\alpha:=\{\alpha_{t}\}_{t \geq 0}$. The associated product system is denoted by $E:=(E(t))_{t>0}$ and the set of decomposable vectors in $E(t)$ is denoted by $D(t)$. Fix $t>0$. For $\xi \in Ker(V_{t}^{*})$, let $T^{(t)}_{e(\xi)} \in B(\Gamma(\clh))$ be defined by
\begin{equation}
\label{exponential}
T^{(t)}_{e(\xi)}e(\eta):=e(\xi+V_{t}\eta)
\end{equation}
for $\eta \in \clh$. For $t>0$, let $e_{t}:=\Gamma(V_{t})$. Then $(e_{t})$ is  a left coherent section of decomposable vectors. 
The following are straightforward consequences of the definitions.
\begin{enumerate}
\item[(1)] For $t>0$ and $\xi \in Ker(V_{t}^{*})$, $T^{(t)}_{e(\xi)} \in D(t)$.
\item[(2)] For $t>0$ and $\xi,\eta \in Ker(V_{t}^{*})$, $L^{e}(t,T^{(t)}_{e(\xi)},T^{(t)}_{e(\eta)})=\langle \xi|\eta \rangle$.
\end{enumerate}

Fix $t>0$. Let $\mathcal{K}(t)$ be the Hilbert space built out of the usual GNS construction applied to the pair $(D(t), L^{e}(t,.,.))$. That is, there exists a map $F:D(t) \to \mathcal{K}(t)$ such that 
\begin{enumerate}
\item[(1)] for $x,y \in D(t)$, $\langle F(x)|F(y) \rangle =L^{e}(t,\dot{x},\dot{y})$, and
\item[(2)] the set $\{F(x):x \in D(t)\}$ is total in $\mathcal{K}(t)$.
\end{enumerate}

Note that the map $Ker(V_{t}^{*}) \ni \xi \to F(T^{(t)}_{e(\xi)}) \in \mathcal{K}(t)$ preserves the inner product. As a consequence, it follows that the map $Ker(V_{t}^{*}) \ni \xi \to F(T^{(t)}_{e(\xi)}) \in \mathcal{K}(t)$ is a linear isometry. With the foregoing notation, we have the following description of decomposable vectors of CCR flows.

\begin{ppsn}
\label{decomposable}
Fix $t>0$. Suppose $x \in D(t)$. Then there exists $\xi \in Ker(V_{t}^{*})$ and $\lambda \in \mathbb{C}\backslash \{0\}$ such that $x=\lambda T^{(t)}_{e(\xi)}$.
\end{ppsn}
\textit{Proof.} It suffices to prove assuming $\langle x|e_{t} \rangle =1$.  By the preceding remarks, it follows that the map $Ker(V_{t}^{*}) \ni \eta \to L^{e}(t,T^{(t)}_{e(\eta)},x) \in \mathbb{C}$ is a bounded linear functional. Thus there exists $\xi \in Ker(V_{t}^{*})$ such that
$L^{e}(t,T^{(t)}_{e(\eta)},x)=\langle \eta|\xi \rangle$ for every $\eta \in Ker(V_{t}^{*})$. Taking exponentials, we obtain $\langle T^{(t)}_{e(\eta)}|x \rangle = \langle T^{(t)}_{e(\eta)}|T^{(t)}_{e(\xi)} \rangle$. Since $\{T^{(t)}_{e(\eta)}: \xi \in Ker(V_{t}^{*})\}$ is total in $E(t)$, it follows that $x=T^{(t)}_{e(\xi)}$. This completes the proof. \hfill $\Box$



\section{Arveson's $e$-logarithm}
In this section, we construct    Arveson's $e$-logarithm in the higher dimensional setting. 
Let $P \subset \mathbb{R}^{d}$ be a closed convex cone which we assume to be spanning and pointed, i.e. $P-P=\mathbb{R}^{d}$ and $P \cap -P=\{0\}$. The cone $P$ will be fixed for the remainder of this paper. Let us fix a few notation that we will use throughout. The letter $\Omega$ stands for the interior of $P$. Note that $\Omega$ is an ideal in $P$, i.e. $\Omega+P \subset \Omega$. Also $\Omega$ is dense in $P$. For $x, y \in \mathbb{R}^{d}$, we write $x \leq y$ if $y-x \in P$ and $x<y$ if $y-x \in \Omega$. We have the following \textbf{Archimedean property:} Let $a \in \Omega$ be given. Then given $x \in \mathbb{R}^{d}$, there exists $n \geq 1$ such that $x<na$.

Let $\alpha:=\{\alpha_{x}\}_{x \in P}$ be an $E_0$-semigroup over $P$ on $B(\clh)$. For $x \in P$, let 
\[
E(x):=\{T \in B(\clh): \alpha_{x}(A)T=TA , \forall A \in B(\clh)\}.
\]
For $x \in P$, $E(x)$ is a separable Hilbert space with the inner product given by $\langle S|T \rangle=T^{*}S$ for $S,T \in E(x)$. Just like in the one dimensional case, the bundle of Hilbert spaces $E:=\{E(x):x \in P\}$ has an associative multiplication and an appropriate measurable structure. We call $E$  the product system associated to $\alpha$. For more details regarding  the exact connection between multiparameter $E_0$-semigroups and product systems,  we refer the reader to \cite{Murugan_Sundar_continuous} and \cite{Murugan_Sundar}.

Just like in the one dimensional case, we define the notion of decomposability as follows. Fix $x \in P$ and $u \in E(x)$. We say that $u$ is a decomposable vector if $u$ is non-zero and given $y \in P$ with $y \leq x$, there exists $v \in E(y)$ and $w \in E(x-y)$ such that 
$u=vw$. We denote the set of decomposable vectors in $E(x)$ by $D(x)$. 

\begin{rmrk}
\label{A few remarks}
\begin{enumerate}
\item[(0)] Fix $x \in P$ and $u \in D(x)$. It is clear that $u$ is a decomposable vector for the $1$-parameter product system $\{E(tx):t >0\}$.
\item[(1)] The proof of Proposition 6.0.2 of \cite{Arveson} carries over in the multidimensional context as well and we have the following. 
Suppose $x \in P$ and $u \in D(x)$. Let $y,z \in P$ be such that $y+z=x$ and let $v \in E(y)$, $w \in E(z)$ be such that $u=vw$. Then $v \in D(y)$ and $w \in D(z)$. 
\item[(2)] However we are unable to settle the following fundamental question. Let $x,y \in P$ and $u\in D(x)$, $v \in D(y)$ be given. Is it true that $uv \in D(x+y)$? The main difficulty is that unlike in the case of the real numbers, the ordering induced by the cone $P$ is not a total order.  The second difficulty is that we are yet to construct enough examples in the multiparameter context to test the hypothesis.  One good example is the following. Using the  CAR construction, we could construct an $E_0$-semigroup starting from an isometric representation on the algebra of bounded operators on the antisymmetric Fock space. To test the hypothesis, it is essential to determine the decomposable vectors of one parameter CAR flows in a coordinate free manner. We consider this as an interesting and an important problem.
\item[(3)] The notions of left( right) divisors, left (right) coherent sections and propagators are defined exactly the same way as in the one dimensional case and we do not repeat the definitions. For example, let $(u_{x})_{x \in P}$ be a section of decomposable vectors. We say that $(u_{x})_{x \in P}$ is left coherent if given $x,y \in P$ with $y \leq x$, there exists a necessarily unique element $u(y,x) \in D(x-y)$, called a \textbf{propagator}, such that $u_{x}=u_{y}u(y,x)$.
\end{enumerate}
\end{rmrk}

Following Arveson, we call an $E_0$-semigroup  \textbf{decomposable} if the following two conditions are satisfied.
\begin{enumerate}
\item[(1)] For $x,y \in P$, $D(x)D(y) \subset D(x+y)$, and
\item[(2)] for every $x \in P$, $D(x)$ is total in $E(x)$.
\end{enumerate}
Condition $(1)$ is equivalent to the assertion that for $x,y \in P$, $D(x)D(y)=D(x+y)$. Note that Condition $(1)$ is automatically satisfied when $P=[0,\infty)$.

We  first prove that CCR flows are indeed decomposable. Let $V:=(V_{x})_{x \in P}$ be an isometric representation on a Hilbert space $\clh$ and denote the corresponding CCR flow on $B(\Gamma(\clh))$ by $\alpha:=\{\alpha_{x}\}_{x \in P}$. For $x \in P$, set $E_{x}:=V_{x}V_{x}^{*}$ and $E_{x}^{\perp}=1-E_{x}$. Fix $x \in P$ and let $\xi \in Ker(V_{x}^{*})$ be given. Then the ``exponential vectors" $T^{(x)}_{e(\xi)}$ on $B(\Gamma(\clh))$ are defined exactly as in Eq. \ref{exponential}. Denote the product system associated to $\alpha$ by $E$. 

Fix $x \in P$ and $\xi \in Ker(V_{x}^{*})$. Note that for $y \in P$ with $y \leq x$, we have the equality $T^{(x)}_{e(\xi)}=T^{(y)}_{e(E_{y}^{\perp}\xi)}T^{(x-y)}_{e(V_{y}^{*}\xi)}$. This implies that $T^{(x)}_{e(\xi)}$ is a decomposable vector. This together with Proposition \ref{decomposable}  and Remark \ref{A few remarks} implies that \[D(x)=\{\lambda T^{(x)}_{e(\xi)}: \lambda \in \mathbb{C} \backslash \{0\}, \xi \in Ker(V_{x}^{*})\}.\] Noting that the product of ``exponential" vectors is again an ``exponential" vector, we obtain the following.

\begin{ppsn}
Let $V$ be an isometric representation on a Hilbert space $\clh$. Then the CCR flow $\alpha^{V}$ associated to $V$ on $B(\Gamma(\clh))$ is decomposable. 

\end{ppsn} 

Fix a decomposable $E_{0}$-semigroup $\alpha$ for the remainder of this section. 
We proceed towards proving the existence of the $e$-logarithm.  We imitate Arveson. The first step in this direction is to address the existence of left coherent sections. For $x \in P$, define an equivalence relation on $D(x)$ by identifying two vectors if they are scalar multiples of each other and denote the set of equivalence classes by $\Delta(x)$. Then $\displaystyle \Delta:=\coprod_{x \in P}\Delta(x)$ is a ``path space" over $P$. For $u \in D(x)$, we denote the equivalence class in $\Delta(x)$ by $\dot{u}$.
\begin{ppsn}
\label{existence of left coherent sections}
Let $a \in \Omega$ and $e \in D(a)$ be of norm one. Then there exists a left coherent section $(e_{x})_{x \in P}$ of decomposable vectors of norm one such that $e_{a}=e$.
\end{ppsn}
\textit{Proof.} For $n \geq 1$, let $E_{n}=\{x \in P: 0 \leq x \leq na \}$. Note that $E_{n}$ is increasing and by the Archimedean property, we have $\bigcup_{n=1}^{\infty}E_{n}=P$.   For every $n$, $\dot{e^{n}} \in \Delta(na)$. Thus for $x \in E_{n}$, we have  unique elements $\dot{e^{(n)}_{x}} \in \Delta(x)$ and $\dot{f^{(n)}_{x}} \in \Delta(na-x)$ such that $\dot{e^{n}}=\dot{e^{(n)}_{x}}\dot{f^{(n)}_{x}}$. It is routine to verify that $\{\dot{e^{(n)}_{x}}: n \geq 1, x \in E_{n}\}$ patches together to define a well defined left coherent section  $(\dot{e_{x}})_{x \in P}$ of $\Delta$ such that $\dot{e_{a}}=\dot{e}$.  For each $x \in P$, choose a representative $e_{x}$ of $\dot{e_{x}}$ of norm one such that $e_{a}=e$. Then $(e_{x})_{x \in P}$ is the desired left coherent section. This completes the proof. \hfill $\Box$

\begin{lmma}
\label{crucial sequence}
Let $(x_n)$ be a sequence in $P$. Suppose $(u_n), (v_n)$ are two sequences  of decomposable vectors of norm one such that $u_n,v_n \in D(x_n)$. Suppose that $(x_n)$ decreases to $0$, i.e. $x_{n+1} \leq x_{n}$ and $x_{n} \to 0$. Assume that $(u_n)$ and $(v_n)$ are left coherent, i.e. there exist $\widetilde{u}_{n}, \widetilde{v}_{n} \in D(x_n-x_{n+1})$ such that $u_{n}=u_{n+1}\widetilde{u}_{n}$ and $v_{n}=v_{n+1}\widetilde{v}_{n}$. Then $|\langle u_{n}|v_{n} \rangle| \to 1$. 
\end{lmma}
\textit{Proof.} Choose $a \in \Omega$ such that $x_n <a$ for every $n$. Pick an element $w \in D(a-x_1)$ and set $e:=u_{1}w \in D(a)$. For $x \leq a$, choose $e_{x} \in D(x)$ of unit norm such that $e_{x}$ is a left divisor of $e$. We can arrange in a such a way that $e_{x_n}=u_{n}$.  Similary, for $x \leq a$, choose $f_{x} \in D(x)$ of unit norm such that $f_{x}$ is a left divisor of $f:=v_{1}w$ and $f_{x_{n}}=v_{n}$. We claim that there exists a subsequence of $(|\langle u_{n}|v_{n} \rangle|)_{n \geq 1}$ which converges to $1$.

For every $k$, $0<\frac{a}{k}$. Hence there exists a subsequence $(x_{n_k})$ such that $x_{n_k}<\frac{a}{k}$. Note that $u_{n_k}$ is a left divisor of $e_{\frac{a}{k}}$ and both have norm $1$. Similarly, $v_{n_k}$ is a left divisor of $f_{\frac{a}{k}}$. Hence 
\[
|\langle e_{\frac{a}{k}}|f_{\frac{a}{k}} \rangle|  \leq |\langle u_{n_k}|v_{n_k} \rangle| \leq 1.
\]
By Theorem 6.1.1 of \cite{Arveson}, we have $|\langle e_{\frac{a}{k}}|f_{\frac{a}{k}}\rangle| \to 1$. Now the above inequality implies that $|\langle u_{n_k}|v_{n_k} \rangle| \to 1$. This proves our claim.

Repeating the above argument for a subsequence for $|\langle u_n|v_n \rangle|$, we conclude that every subsequence of $|\langle u_n|v_n \rangle|$ has a further subsequence which converges to $1$. Hence $|\langle u_n|v_n \rangle| \to 1$. The proof is now complete. \hfill $\Box$.

\begin{lmma}
\label{decreasing in the interior}
Let $(x_n)$ be a sequence in $\Omega$ such that $x_{n}$ converges to $0$. Then there exists a subsequence $(x_{n_k})$ which decreases to $0$.
\end{lmma}
\textit{Proof.} Choose $a \in \Omega$. Since $0 <a$, there exists $x_{n_1}$ such that $x_{n_1}<a$. Define $x_{n_k}$ inductively as follows. Note that $0<x_{n_{k-1}}$ and $0<\frac{a}{k}$. Hence there exists $x_{n_k}$ such that $x_{n_k}<x_{n_{k-1}}$ and $x_{n_k} <\frac{a}{k}$. It is now clear that $(x_{n_k})$ decreases to $0$. This completes the proof. \hfill $\Box$

Let us recall the notion of a dual cone. For details, we refer the reader to  \cite{Faraut}.  The dual  of $P$ denoted $P^{*}$ is defined by  
\[
P^{*}:=\{y \in \mathbb{R}^{d}: \langle x|y \rangle \geq 0,~ \forall x \in P\}.
\]
It is well known that $P^{*}$ is pointed and spanning. Moreover the dual of $P^{*}$ is $P$. 
 The cone $P$ is said to be polyhedral if there exists $v_1,v_2,\cdots,v_k \in \mathbb{R}^{d}$ such that \[P=\{\sum_{i=1}^{k}r_{i}v_{i}: r_{i} \geq 0\}.\] Farkas' theorem states that if $P$ is polyhedral,  then the dual $P^{*}$ is also polyhedral (See Page 11 of \cite{Fulton}).

\begin{lmma}
\label{decreasing}
Suppose that $P$ is polyhedral. Let $(z_n)$ be a sequence in $P$ which conveges to $0$. Then there exists a subsequence $(z_{n_k})$ which decreases to $0$.
\end{lmma}
\textit{Proof.} Let $v_{1},v_{2}, \cdots, v_{\ell} \in P^{*}$ be such that $P^{*}=\{\sum_{i=1}^{\ell}r_{i}v_{i}: r_{i} \geq 0\}$. Choose a subsequence $(z_{n_k})$ such that for every $i=1,2,\cdots, \ell$, $\langle z_{n_k}|v_{i} \rangle$ decreases to $0$. Let $v \in P^{*}$ be given. Write $v=\sum_{i=1}^{\ell}r_{i}v_{i}$ with $r_i \geq 0$. It is clear that $\langle z_{n_k}|v \rangle$ decreases to $0$. This implies that for every $k$ and $v \in P^{*}$, $\langle z_{n_k}-z_{n_{k+1}}|v \rangle \geq 0$. Since the dual of $P^{*}$ is $P$, it follows that $z_{n_{k+1}}-z_{n_k} \in P$. Consequently, the sequence $(z_{n_k})$ decreases to $0$. This completes the proof. \hfill $\Box$

\begin{rmrk}
Lemma \ref{decreasing} is not true in general as the following example shows. Let $V$ be a finite dimensional real inner product space of dimension at least $2$. Denote the space of symmetric operators on $V$ by $\mathcal{S}(V)$ and the cone of positive operators on $V$ by $\mathcal{P}(V)$. Then $\mathcal{P}(V)$ is indeed a closed spanning cone in $\mathcal{S}(V)$. Moreover $\mathcal{P}(V)$ is pointed. Choose a sequence $(v_n) \in V$ of unit vectors such that $v_{m}$ is not a scalar multiple of $v_{n}$ if $ m \neq n $. Denote the orthogonal projection onto $span\{v_n\}$ by $P_{n}$. Set $z_{n}=\frac{1}{n}P_n$. Note that $z_n \to 0$. It is not difficult to see that no subsequence of $z_n$ is decreasing.
\end{rmrk}

\begin{ppsn}
\label{crucial}
Let $(x_n)$ be a sequence in $P$  and $a \in \Omega$. Suppose that $(x_n) \to a $. Then there exists a subsequence $(x_{n_k})$ and a sequence $(t_{n_k})$ of positive real numbers such that the following holds.
\begin{enumerate}
\item[(1)] The sequence $(t_{n_k}) \to 1$.
\item[(2)] If we set $z_{n_k}=x_{n_k}-t_{n_k}a$ then $z_{n_k} \in P$ and decreases to $0$.
\end{enumerate}
\end{ppsn}
\textit{Proof.}  With no loss of generality, we can assume that $x_{n} \in \Omega$. 

\textbf{Step 1:} First we prove the statement assuming that $P$ is polyhedral. Denote the boundary of $P$ by $\partial(P)$. Prop. 2.3 of \cite{Anbu_Sundar} implies that the map \[\partial(P) \times (0,\infty) \ni (z,t) \to z+ta \in \Omega\] is a homeomorphism. Hence there exists $(z_n) \in \partial(P)$ and $t_{n} \in (0,\infty)$ such that $z_n \to 0$ and $t_n \to 1$ with $x_n=z_n+t_na$. The desired conclusion is immediate if we apply Lemma \ref{decreasing} to the sequence $(z_n)$. 

\textbf{Step 2:} Next we consider the situation when $P$ is not necessarily polyhedral. Choose a basis $v_{1},v_{2},\cdots,v_{d}$ of $\mathbb{R}^{d}$ and set $v_{d+1}:=-\sum_{i=1}^{d}v_{i}$. Let $k$ be a large natural number such that for every $i=1,2,\cdots,d+1$, $a+\frac{1}{k}v_{i} \in \Omega$. This is possible as $\Omega$ is open and $a \in \Omega$. For $i=1,2,\cdots, d+1$, set $\widetilde{v_i}=a+\frac{1}{k}v_{i}$.  Define \[
Q:=\{\sum_{i=1}^{d+1}r_{i}\widetilde{v}_{i}: r_{i} \geq 0\}.\] 
Note that $\{\widetilde{v}_{i}:i=1,2\cdots,d+1\}$ spans $\mathbb{R}^{d}$. Thus $Q$ is spanning. Moreover $Q$ is pointed as $Q$ is contained in $P$. Observe that the equality $a:=\sum_{i=1}^{d+1}\frac{1}{d+1}\widetilde{v}_{i}$ implies that $a$ is in the interior of $Q$. Eventually, $x_{n}$ lies in the interior of $Q$. The desired conclusion follows by applying Step 1 to the sequence $(x_n)$ lying in the polyhedral cone $Q$. 
 This completes the proof. \hfill $\Box$

The main theorem that allows us to construct the $e$-logarithm is the following.
\begin{thm}
\label{main theorem}
Let $(u_{x})_{x \in P}$ and $(v_{x})_{x \in P}$ be two left coherent sections of decomposable vectors with unit norm. Then the map $\Omega \cup \{0\} \ni x \to |\langle u_{x}|v_{x} \rangle| \in \mathbb{C}$ is continuous.
\end{thm}
\textit{Proof.} Denote the propagators of $(u_{x})_{x \in P}$ by $\{u(y,x):0 \leq y \leq x\}$ and the propagators of $(v_{x})_{x \in P}$ by $\{v(y,x):0 \leq y \leq x\}$.

 First we consider continuity at an interior point. Fix $a \in \Omega$ and let $(x_{n})$ be a sequence in $\Omega$ such that $(x_{n}) \to a$. We  show that there exists a subsequence of $(|\langle u_{x_n}|v_{x_n} \rangle|)$ which converges to $|\langle u_{a}|v_{a} \rangle|$. 
Choose a subsequence $(x_{n_k})$ and a sequence $(t_{n_k})$ as in Prop. \ref{crucial}. Set $\widetilde{x}_{k}:=x_{n_k}$, $s_{k}:=t_{n_k}$ and $y_k:=\widetilde{x}_{k}-s_{k}a$. Then $(y_k)$ decreases to $0$ and $s_{k} \to 1$.

For $k=1,2,\cdots$, let $u_{k}(s)=u(y_{k},y_{k}+sa)$ and $v_{k}(s)=v(y_{k},y_{k}+sa)$. Then $u_{k}$ and $v_{k}$ are left coherent sections of decomposable vectors with unit norm in the $1$-parameter product system $\{E(sa):s \geq 0\}$. Define $f_{k}:[0,\infty) \to \mathbb{C}$ by \[
f_{k}(t)=| \langle u_{k}(t)|v_{k}(t) \rangle |= |\langle u(y_{k},y_{k}+ta)|v(y_{k},y_{k}+ta)\rangle|.
\]
Theorem 6.1.3 of \cite{Arveson} asserts that for every $k$, $f_{k}$ is continuous. Let $f:[0,\infty) \to \mathbb{C}$ be defined by $f(t)=|\langle u(0,ta)|v(0,ta) \rangle|$. A similar reasoning implies that $f$ is continuous. 

We claim that $f_{k}$ converges pointwise to $f$. Fix $t \geq 0$. Observe that
\begin{align*}
u(0,ta)u(ta,y_{k}+ta)&=u(0,y_{k}+ta) \\
                                &=u(0,y_{k})u(y_{k},y_{k}+ta).
\end{align*}
A similar equation holds for $v$. Thus we have
\begin{equation}
\label{convergence}
|\langle u_{ta}|v_{ta}\rangle||\langle u(ta,y_{k}+ta)|v(ta,y_{k}+ta)\rangle|=|\langle u_{y_{k}}|v_{y_{k}}\rangle||\langle u(y_{k},y_{k}+ta)|v(y_{k},y_{k}+ta)\rangle|.
\end{equation}

By  Lemma \ref{crucial sequence},  it follows that $|\langle u(ta,y_{k}+ta)|v(ta,y_{k}+ta)\rangle| \to 1$ and $|\langle u_{y_k}|v_{y_k} \rangle| \to 1$. Now Eq. \ref{convergence} implies that $f_{k}(t) \to f(t)$. This proves our claim. 

Let $\epsilon_{k}=|\langle u(y_{k+1},y_{k})|v(y_{k+1},y_{k})\rangle|$. By Cauchy-Schwarz inequality, it follows that $\epsilon_{k} \leq 1$. Note that $\epsilon_{k}=\frac{|\langle u_{y_k}|v_{y_k}\rangle|}{|\langle u_{y_{k+1}}|v_{y_{k+1}}\rangle|}$. By Lemma \ref{crucial sequence}, it follows that $\epsilon_{k} \to 1$. By passing to a subsequence, if necessary, we can assume that $|\log(\epsilon_k)|<\frac{1}{2^{k}}$. Define $\widetilde{\epsilon}_{k}=e^{-\frac{1}{k}}$. Note that there exists $N$ such that for $k \geq N$,
\begin{equation}
\label{estimation}
\frac{\widetilde{\epsilon}_{k}}{\widetilde{\epsilon}_{k+1}} \leq \epsilon_{k}.
\end{equation}
For $k \geq 1$, set $g_{k}:=\widetilde{\epsilon}_{k}f_{k}$. We claim that for every $t \in [0,\infty)$, $(g_{k}(t))_{k=N}^{\infty}$ is an increasing sequence. Fix $t \geq 0$. Calculate as follows to observe that
\begin{align*}
\epsilon_{k}f_{k}(t)&=|\langle u(y_{k+1},y_{k})|v(y_{k+1},y_{k})\rangle||\langle u(y_{k},y_{k}+ta)|v(y_{k},y_{k}+ta)\rangle| \\
                             &=|\langle u(y_{k+1},y_{k}+ta)|v(y_{k+1},y_{k}+ta) \rangle| \\
                             &=|\langle u(y_{k+1},y_{k+1}+ta)u(y_{k+1}+ta,y_{k}+ta)|v(y_{k+1},y_{k+1}+ta)v(y_{k+1}+ta,y_{k}+ta)\rangle| \\
                             &=|\langle u(y_{k+1},y_{k+1}+ta)|v(y_{k+1},y_{k+1}+ta)\rangle||\langle u(y_{k+1}+ta,y_{k}+ta)|v(y_{k+1}+ta,y_{k}+ta)\rangle |\\
                             & \leq f_{k+1}(t)||u(y_{k+1}+ta,y_{k}+ta)||||v(y_{k+1}+ta,y_{k}+ta)||\\
                             &=f_{k+1}(t).
\end{align*}
The above calculation together with Eq. \ref{estimation} implies that for $k \geq N$, $g_{k}(t) \leq g_{k+1}(t)$ for every $t \in [0,\infty)$. Also observe that since $\widetilde{\epsilon}_{k} \to 1$, it follows that $g_{k}$ converges pointwise to $f$. 
By Dini's theorem, we conclude that $g_{k}$ converges to $f$ uniformly on compact subsets of $[0,\infty)$. As a consequence, we obtain that $f_{k}$ converges uniformly to $f$ on compact subsets of $[0,\infty)$.

Calculate as follows to observe that 
\begin{equation}
\label{continuity}
|\langle u_{x_{n_k}}|v_{x_{n_k}}\rangle|=|\langle u(0,y_{k})u(y_k,y_{k}+s_{k}a)|v(0,y_{k})v(y_{k},y_{k}+s_{k}a)\rangle| \\
=|\langle u_{y_k}|v_{y_k} \rangle|f_{k}(s_k).
\end{equation}
Lemma \ref{crucial sequence} implies that $|\langle u_{y_k}|v_{y_k}\rangle| \to 1$. The uniform convergence of $f_{k}$ to $f$ on compact sets implies that $f_{k}(s_k) \to f(1)=|\langle u_{a}|v_{a}\rangle|$. From Eq. \ref{continuity}, we have that the sequence 
$|\langle u_{x_{n_k}}|v_{x_{n_k}}\rangle| \to |\langle u_{a}|v_{a} \rangle|$. By repeating the argument for a subsequence of $(|\langle u_{x_n}|v_{x_n} \rangle|)$, we conclude that every subsequence of $(|\langle u_{x_n}|v_{x_n} \rangle|)$ has a further subsequence which converges to $|\langle u_{a}|v_{a} \rangle|$. Hence $(|\langle u_{x_n}|v_{x_n} \rangle|) \to |\langle u_a|v_a \rangle|)$. 

Continuity at the origin is easier. Let $(x_n)$ be a sequence in $\Omega$ such that $(x_n) \to 0$. Use Lemma \ref{decreasing in the interior} to  choose a subsequence $(x_{n_k})$ which decreases to $0$. Then Lemma \ref{crucial sequence} implies that $|\langle u_{x_{n_k}}|v_{x_{n_k}}\rangle| \to 1$. Thus $(|\langle u_{x_n}|v_{x_n} \rangle|)$ has a subsequence which converges to $1$. By our usual argument, i.e. by repeating the above argument for a subsequence of $(|\langle u_{x_n}|v_{x_n}\rangle|)$, we conclude that every subsequence of $(|\langle u_{x_n}|v_{x_n}\rangle|)$ has a further subsequence which converges to $1$. Consequently $(|\langle u_{x_n}|v_{x_n}\rangle|) \to 1$. This completes the proof. \hfill $\Box$

Let $D$ denote the set of left coherent sections of decomposable vectors with unit norm. Consider an element $e=(e_{x})_{x \in P} \in D$. Denote by $D^{e}$ the set of left coherent sections of decomposable vectors $(u_{x})_{x \in P}$ satisfying the equation $\langle u_x|e_x \rangle=1$ for every $x \in P$. 

\begin{lmma}
\label{continuity of inner product}
Let $u,v \in D^{e}$ be given. Then we have the following.
\begin{enumerate}
\item[(1)] The map $\Omega \cup \{0\} \ni  x \to ||u_{x}|| \in (0,\infty)$ is continuous and increasing with respect to the order induced by $P$.  Moreover $||u_{x}|| \geq 1$ for every $x \in \Omega \cup \{0\}$.
\item[(2)] The map $\Omega \cup \{0\} \ni x \to \langle u_{x}|v_{x} \rangle \in \mathbb{C}$ is continuous.
\end{enumerate}
\end{lmma}
\textit{Proof.} The conclusion of $(1)$ follows by applying the argument outlined in Lemma 6.3.3 of \cite{Arveson} and by making use of Theorem \ref{main theorem}.

Let $(x_n)$ be a sequence in $\Omega$ such that $(x_n) \to a \in \Omega$. Without loss of generality, we can assume that $x_{n}<2a$. We show that there exists a subsequence of $(\langle u_{x_n}|v_{x_n}\rangle)$ which converges to $\langle u_{a}|v_{a} \rangle$. Choose a sequence of real numbers $(t_{n})$ such that $2>t_{n}>1$ and $t_{n} \to 1$. For every $k$, $a<t_{k}a$. Thus there exists a subsequence $(x_{n_k})$ such that $x_{n_k}<t_{k}a$. Arguing as in Theorem 6.3.4 of \cite{Arveson}, we obtain the following estimate
\begin{equation}
\label{estimate}
|\langle u_{x_{n_k}}|v_{x_{n_k}}\rangle - \langle u_{t_ka}|v_{t_ka} \rangle| \leq ||u_{2a}||||v_{2a}||\sqrt{(||u_{t_ka}||^{2}-||u_{x_{n_k}}||^{2})(||v_{t_ka}||^{2}-||v_{x_{n_k}}||^{2})}.
\end{equation}
Eq. \ref{estimate} and Part (1) of this Lemma implies that $\langle u_{x_{n_k}}|v_{x_{n_k}}\rangle-\langle u_{t_ka}|v_{t_ka}\rangle \to 0$. Theorem 6.3.4 of \cite{Arveson} implies that $\langle u_{t_ka}|v_{t_ka} \rangle \to \langle u_{a}|v_{a} \rangle$. Hence a subsequence of $(\langle u_{x_n}|v_{x_n} \rangle)$ converges to $\langle u_{a}|v_{a} \rangle$. Repeating the same argument again for a subsequence of $(\langle u_{x_n}|v_{x_n} \rangle)$, we deduce that every subsequence of $(\langle u_{x_n}|v_{x_n} \rangle)$ has a further subsequence which converges to $\langle u_{a}|v_{a} \rangle$. This implies that $\langle u_{x_n}|v_{x_n} \rangle \to \langle u_{a}|v_{a} \rangle$. 

Continuity at the origin follows from the following estimate. For $x \in \Omega$, 
\begin{align*}
|\langle u_{x}|v_{x} \rangle -1|&=|\langle u_{x}-e_{x}|v_{x}-e_{x} \rangle| \\
                                              & \leq ||u_{x}-e_{x}||||v_{x}-e_{x}||\\
                                              & \leq \sqrt{||u_{x}||^{2}-1}\sqrt{||v_{x}||^{2}-1}.
\end{align*}
This completes the proof. \hfill $\Box$

With Lemma \ref{continuity of inner product} in hand, we can  define the $e$-logarithm. Restrict the path space $\Delta$ which is apriori defined over $P$ to the subsemigroup $\Omega$ which we again denote by $\Delta$. 
Let \[\Delta^{(2)}:=\{(x,\dot{u},\dot{v}):x \in \Omega, u,v \in D(x)\}.\]
Let us recall Arveson's notion of continuity for functions defined on $\Delta$. Consider a function $\phi: \Delta \to \mathbb{C}$. We say that $\phi$ is continuous if for every left coherent section $(u_{x})_{x \in P}$ of decomposable vectors, the map $\Omega \ni x\to \phi(x,\dot{u}_x) \in \mathbb{C}$ is continuous. We say that $\phi$ vanishes at the origin if for every left coherent section $(u_{x})_{x \in P}$, $\lim_{x \to 0}\phi(x,\dot{u}_x) =0$.  We make similar definitions for functions defined on $\Delta^{(2)}$.

\begin{thm}
\label{e-logarithm}
Let $e=(e_{x})_{x \in P} \in D$ be given. Then there exists a unique continuous function $L^{e}:\Delta^{(2)} \to \mathbb{C}$ which vanishes at the origin and for $x \in \Omega$ and $u,v \in D(x)$,
\[
e^{L^{e}(x,\dot{u},\dot{v})}=\frac{\langle u| v \rangle}{\langle u|e_{x} \rangle \langle e_{x}|v\rangle }.
\]
Suppose $f=(f_{x})_{x \in P}$ is another element in $D$. Denote the corresponding continuous function defined on $\Delta^{(2)}$ by $L^{f}$. Then there exists a continuous function $\phi: \Delta \to \mathbb{C}$ that vanishes at the origin such that 
\begin{equation}
\label{independence}
L^{f}(x,\dot{u},\dot{v})=L^{e}(x,\dot{u},\dot{v})+\phi(x,\dot{u})+\overline{\phi(x,\dot{v})}
\end{equation}
for $x \in \Omega$ and $u,v \in D(x)$. 
\end{thm}
\textit{Proof.} We merely give a sketch of the proof as the proof is essentially the same as that of Theorem 6.4.2 of \cite{Arveson}. Fix $a \in \Omega$ and $u,v \in D(a)$. By Prop. \ref{existence of left coherent sections}, there exist left coherent sections $(u_{x})_{x \in P}$ and $(v_{x})_{x \in P}$ of decomposable vectors such that $u_{a}=u$ and $v_{a}=v$. By Lemma \ref{continuity of inner product}, it follows that the map 
\[
\Omega \cup \{0\} \ni x \to \frac{\langle u_{x}|v_{x} \rangle}{\langle u_{x}|e_{x} \rangle \langle e_{x}|v_{x} \rangle} \in \mathbb{C}\backslash \{0\}\] is continuous and takes value $1$ at $x=0$. 
Note that $\Omega \cup \{0\}$ is a contractible topological space. Thus there exists a unique continuous function $l: \Omega \cup \{0\} \to \mathbb{C}$ such that $l(0)=0$ and 
\[
e^{l(x)}= \frac{\langle u_{x}|v_{x} \rangle}{\langle u_{x}|e_{x} \rangle \langle e_{x}|v_{x} \rangle}.\]
Set $L^{e}(a,\dot{u},\dot{v}):=l(a)$.  The well-definedness of $L^{e}$ and the other conclusions follow as in Theorem 6.4.2 of \cite{Arveson} by making repeated use of Theorem \ref{main theorem} and Lemma \ref{continuity of inner product}.  \hfill $\Box$

\begin{rmrk}
\label{comparision}
Fix $e=(e_{x})_{x \in P} \in D$. Let $a \in \Omega$ be given. Note that $\widetilde{e}:=\{e_{ta}\}_{t>0}$ is left coherent section of the one parameter product system $\{E(ta):t>0\}$. Let $L^{\widetilde{e}}$ be the $e$-logarithm corresponding to $\widetilde{e}$ associated to the $1$-parameter decomposable product system $\{E(ta):t>0\}$. It is clear from the definitions that for $u,v \in D(a)$, 
\[L^{e}(a,u,v)=L^{\widetilde{e}}(1,u,v).\]
An immediate consequence of the above equation and  Theorem 6.5.1. of \cite{Arveson} is that for every $a \in \Omega$, the map \[D(a) \times D(a) \to L^{e}(a,\dot{u},\dot{v}) \in \mathbb{C}\] is positive definite. 
\end{rmrk}

Fix $e=(e_{x})_{x \in P} \in D$ and let $L^{e}$ be the corresponding $e$-logarithm. There is a small subtlety involved in proving the additivity of $L^{e}$

\begin{ppsn}
\label{additivity}
Fix $a \in \Omega$. Then there exists a continuous function $\psi_{a}:\Delta \to \mathbb{C}$ which vanishes at the origin such that 
for $u_1,u_2 \in D(a)$, $b \in \Omega$ and $v_1,v_2 \in D(b)$, 
\[
L^{e}(a+b,u_1v_1,u_2v_2)-L^{e}(a,u_1,u_2)-L^{e}(b,v_1,v_2)=\psi_{a}(b,v_1)+\overline{\psi_{a}(b,v_2)}.\]
\end{ppsn}
\textit{Proof.} Let $\{e(y,x): 0 \leq y \leq x\}$ be the propagators of $e$. By making an appeal to Theorem \ref{main theorem} and Lemma \ref{continuity of inner product}, construct a continuous function $\psi_{a}$ that vanishes at the origin, as in Prop. 6.4.5 of \cite{Arveson}, such that
\[
e^{\psi_{a}(x,\dot{u})}=\frac{|\langle e(a,a+x)|e_{x}\rangle| \langle u|e_{x} \rangle}{\langle u|e(a,a+x) \rangle \langle e(a,a+x)|e_{x} \rangle}\]
for $x \in \Omega$ and $u \in D(x)$.

Fix $u_1,u_2 \in D(a)$. Let $b \in \Omega$ and $v_1,v_2 \in D(b)$ be given. Choose left  coherent sections $(v^{(1)}_{x})_{x \in P}$ and $(v^{(2)}_{x})_{x \in P}$ such that $v^{(1)}_{b}=v_1$ and $v^{(2)}_{b}=v_2$. For $0 \leq t \leq 1$, let 
\begin{align*}
L(t)&=L^{e}(a+tb,u_1v^{(1)}_{tb},u_2v^{(2)}_{tb})-L^{e}(a,u_1,u_2)-L^{e}(tb,v^{(1)}_{tb},v^{(2)}_{tb}) \\
R(t)&=\psi_{a}(tb,v^{(1)}_{tb})+\overline{\psi_{a}(tb,v^{(2)}_{tb})}.
\end{align*}
 We claim that $L$ is continuous on $[0,1]$ and vanishes at the origin. 

Choose left coherent sections $w^{(1)}$ and $w^{(2)}$ such that $w^{(1)}_{a+b}=u_1v^{(1)}_{b}$ and $w^{(2)}_{a+b}=u_2v^{(2)}_{b}$. \footnote{This is where the subtle point lies. The author does not know how to construct explicitly such left coherent sections which is why one needs to include the axiom $D(x)D(y) \subset D(x+y)$ in the definition of a decomposable $E_0$-semigroup. The left coherent sections $w^{(i)}$ are guaranteed by the defining axiom  $D(x)D(y) \subset D(x+y)$.}
Note that for for $i=1,2$ and $0 \leq x \leq b$, $w^{(i)}_{a+x}$ and $u_{i}v^{(i)}_{x}$ are scalar multiples of each other. Hence for $0 \leq t  \leq 1$, 
\[
L(t)=L^{e}(a+tb,w^{(1)}_{a+tb},w^{(2)}_{a+tb})-L^{e}(a,u_1,u_2)-L^{e}(tb,v^{(1)}_{tb},v^{(2)}_{tb}).\]
It is now clear that $L$ is continuous on $[0,1]$ and vanishes at $0$. A routine calculation reveals that $e^{L(t)}=e^{R(t)}$. Also $R$ is continuous on $[0,1]$  and vanishes at $0$. Hence $L(1)=R(1)$ and the desired conclusion follows. The proof is now complete. \hfill $\Box$

\section{Characterisation of CCR flows}
In this section, we discuss Arveson's characterisation of CCR flows in the multiparameter context. With the $e$-logarithm in hand, Arveson's arguments in the one parameter setting work equally well in the multiparameter context and establish Theorem \ref{characterisation of CCR flows} below. 
 
First, let  us recall Arveson's construction of  an isometric representation associated to a decomposable $E_0$-semigroup. Fix a decomposable $E_0$-semigroup $\alpha=\{\alpha_x\}_{x \in P}$ and denote its product system by $E:=\{E(x)\}_{x \in P}$. Let $\Delta$ be the associated path space over $\Omega$. 
Let $e:=(e_x)_{x \in P}$ be a left coherent section of decomposable vectors of unit norm. The $e$-logarithm associated to the section $e$ is denoted by $L^{e}$. 

\begin{enumerate}
\item[(1)] For $a \in \Omega$, let $H_{a}$ be the Hilbert space obtained as follows: Let $\mathbb{C}_{0}\Delta(a)$ be the set of all finitely supported functions on $\Delta(a)$ whose sum is zero. Define a positive semi-definite inner product on $\mathbb{C}_{0}\Delta(a)$ by the following formula:
\[
\langle f|g \rangle = \sum_{u,v \in \Delta(a)}f(u)\overline{g(v)}L^{e}(a,u,v)
\]
for $f,g \in \mathbb{C}_0\Delta(a)$. Note that in view of Eq. \ref{independence}, the sesquilinear form defined above is independent of the chosen section $e$. Then the Hilbert space $H_{a}$ is obtained by completing the genuine  inner product space  that results from $\mathbb{C}_0\Delta(a)$ after passing to the quotient of $\mathbb{C}_0\Delta(a)$ by the subspace of null vectors. 
For $u \in \Delta(a)$, let $\delta_u$ be the characteristic function at $u$.  For $u,v \in \Delta(a)$, the difference $\delta_{u}-\delta_v \in H_a$ is denoted by $[u]-[v]$.
\item[(2)] Let $a,b \in \Omega$ be such that $a<b$. Choose $f \in D(b-a)$. The map \[H_{a} \ni [u]-[v] \to [uf]-[vf] \in H_{b}\] defines an isometry which is independent of the chosen element $f$. We denote this isometry by $V(b,a)$. Let $H_{\infty}$ be the inductive limit of the Hilbert spaces $H_{a}$ where we embedd $H_{a}$ inside $H_b$ via $V(b,a)$ if $a<b$. For $u,v \in D(a)$, the image of $[u]-[v]$ in $H_{\infty}$ under the natural embedding of $H_a$ in $H_{\infty}$ will again be denoted by $[u]-[v]$. 
\item[(3)] Fix $a \in \Omega$. Choose $f \in D(a)$. The maps $H_{b} \ni [u]-[v] \to [fu]-[fv] \in H_{a+b}$ patch up and induce an isometry on $H_{\infty}$. Moreover the isometry thus obtained is independent of the chosen element $f$. We denote this isometry by $V_{a}$. Then $(V_{a})_{a \in \Omega}$ forms a semigroup of isometries. 
\end{enumerate}
 It is clear that the semigroup of isometries $(V_a)_{a \in \Omega}$, up to unitary equivalence, depends only on the cocycle conjugacy class of $\alpha$. In other words, the isometric representation $V=(V_a)_{a \in \Omega}$ is a cocycle conjugacy invariant. 

\begin{ppsn}
\label{separability}
With the foregoing notation, we have the following.
\begin{enumerate}
\item[(1)] The Hilbert space $H_{\infty}$ is separable.
\item[(2)] The semigroup of isometries $(V_a)_{a \in \Omega}$ is strongly continuous.
\item[(3)] The representation $V:=(V_a)_{a \in \Omega}$ is pure i.e. $\bigcap_{a \in \Omega}V_{a}\clh=\{0\}$.

\end{enumerate}

\end{ppsn}
\textit{Proof.} Fix $a \in \Omega$. 
Let $K_{\infty}$ be the Hilbert space obtained via the above process applied to the one parameter product system $\{E(ta):t>0\}$ and $(K_t)_{t>0}$ be the constituent Hilbert spaces.  Denote the $1$-parameter isometric representation obtained for the product system $\{E(ta):t>0\}$ by $(W_t)_{t>0}$. By Remark \ref{comparision}, it follows that 
for every $n \geq 1$, the map
\[
H_{na} \ni [u]-[v] \to [u]-[v] \in K_{n} 
\]
is an isometry. Hence $H_{na}$ is separable for every $n \geq 1$. The Archimedean property states that $\{na: n \geq 1\}$ is cofinal in $\Omega$. Hence  $H_{\infty}$ is the inductive limit of the Hilbert spaces $(H_{na})_{n=1}^{\infty}$. 
This implies that $H_{\infty}$ is separable. 

Moreover, we can view $H_{\infty}$ as a subspace of $K_{\infty}$ and the restriction of $W_{t}$ to $H_{\infty}$ is $V_{ta}$. From the one parameter assertion, it follows that $\displaystyle \bigcap_{t>0}V_{ta}H_{\infty}=\{0\}$. But note that by the Archimedean property, we have $\displaystyle \bigcap_{b \in \Omega}V_{b}H_{\infty}=\bigcap_{t>0}V_{ta}H_{\infty}=0$. Hence the representation $(V_{b})_{b \in \Omega}$ is pure. 

In view of Theorem 10.8.1 of \cite{Hille}, it suffices to verify the strong continuity of the map  $(0,\infty) \ni t \to V_{ta}$. This follows as above from the one parameter assertion.  This completes the proof. \hfill $\Box$

We provide a proof of the fact that isometric representations which are strongly continuous along each ray is strongly continuous in the following proposition which is based on a simple dilation argument.\begin{ppsn}
\label{dilation}
Let $V:=(V_x)_{x \in \Omega}$ be a semigroup of isometries on a separable Hilbert space $\clh$. Suppose that for $a \in \Omega$ and $\xi \in \clh$, the map 
\[
(0,\infty) \ni t \to V_{ta}\xi \in \clh
\]
is continuous. Then for $\xi \in \clh$, the map $\Omega \ni x \to V_{x}\xi \in \clh$ is continuous.\footnote{ The reason for the inclusion of this proposition  
  is because I am unable to understand the proof outlined in \cite{Hille}.}

\end{ppsn}
\textit{Proof.} \textbf{Step 1:} First consider the case of a unitary representation. Let $U:=(U_{x})_{x \in \Omega}$ be a semigroup of unitaries on $\clh$ satisfying the hypothesis of the proposition. Let $z \in \mathbb{R}^{d}$ be given. Write $z=x-y$ with $x,y \in \Omega$. Set $\widetilde{U}_{z}=U_{x}U_{y}^{*}$. It is straightforward to verify that $\widetilde{U}_{z}$ is well defined and $(\widetilde{U}_z)_{z \in \mathbb{R}^{d}}$ is a group of unitaries. We leave it to the reader to convice himself/herself that for $z \in \mathbb{R}^{d}$ and $\xi \in \clh$, the map 
\[
\mathbb{R} \ni t \to \widetilde{U}_{tz}\xi \in \clh\]
is continuous. Using the fact that $\mathbb{R}^{d}$ has a finite basis, it is routine to prove that for $\xi \in \clh$, the map $\Omega \ni x \to \widetilde{U}_{x}\xi \in \clh$ is continuous. 

\textbf{Step 2:} Let $V=(V_{x})$ be a semigroup of isometries satisfying the hypothesis of the proposition. Let $(U,\clk)$ be the ``minimal" unitary dilation \footnote{Such a dilation is guaranteed by an inductive limit construction.} of $V$, i.e. $\clk$ is a Hilbert space containing $\clh$ as a closed subspace such that 
\begin{enumerate}
\item[(1)] $\{U_{a}\}_{a \in \Omega}$ is a semigroup of unitaries on $\Omega$,
\item[(2)] the union $\displaystyle \bigcup_{a \in \Omega}U_{a}^{*}\clh=\clk$ ( this is the minimality condition), and
\item[(3)] for $\xi \in \clh$ and $a \in \Omega$, $U_{a}\xi=V_{a}\xi$.
\end{enumerate}
From the minimality condition, the hypothesis of the proposition and the fact that addition on $\Omega$ is abelian, it follows that for $\xi \in \clh$, $a \in \Omega$, the map 
\[
(0,\infty) \ni t \to U_{ta} \xi
\]
is continuous. Now (3) together with Step 1 yields the desired conclusion. This completes the proof. \hfill $\Box$

\begin{enumerate}
  \item[(1)] In Arveson's proof of the fact that the decomposable product systems are CCR flows, the first step is the construction of the isometric representation explained above.  This step for the case of a closed convex cone is now  achieved.
  \item[(2)] Then he solves several cohomological problems and exhibits an isomorphism between the product system that one started with and the product system associated to the CCR flow corresponding to the obtained isometric representation. The cohomological issues arise because one does not know, to begin with, whether a decomposable product system admits a unit or not. That units exist, in the one parameter setting, and also in abundance is  a pleasant consequence of Arveson's theorem.
  \item[(3)] But if the decomposable product system has a unit to start with, the above mentioned cohomological issues subside and  Arveson's proof works in the higher dimensional setting too. 
 \end{enumerate}
 We explain $(3)$ below. It is appropriate at this point to recall the definition of a unit of an $E_0$-semigroup. 

\begin{dfn}
Let $\alpha:=\{\alpha_{x}\}_{x \in P}$ be an $E_0$-semigroup over $P$ on $B(\clh)$.  A family $u:=\{u_{x}\}_{x \in P}$ of bounded operators on $\clh$ is called a unit of $\alpha$ if 
\begin{enumerate}
\item[(1)] for $x \in P$, $u_x \neq 0$, 
\item[(2)] for $x, y\in P$, $u_xu_y=u_{x+y}$, 
\item[(3)] for $x \in P$ and $A \in B(\clh)$, $\alpha_{x}(A)u_x=u_xA$, and
\item[(4)] for $\xi,\eta \in \clh$, the map $P \ni x \to \langle u_x\xi|\eta \rangle \in \mathbb{C}$ is measurable.
\end{enumerate}
In other words, if we denote  the product system of $\alpha$ by $E:=\{E(x)\}_{x \in P}$ then a unit is a nowhere vanishing measurable cross section of $E$ which is multiplicative. 
\end{dfn}

Let $\alpha:=\{\alpha_{x}\}_{x \in P}$ be a decomposable $E_0$-semigroup and $E$ be the product system of $\alpha$. Suppose $\alpha$ has a unit.  Fix a unit $(e_x)_{x \in P}$. 
We can assume that $||e_x||=1$. Then $(e_x)$ is clearly a left coherent section of decomposable vectors of unit norm. Let $V$ be the isometric representation associated to $E$ and denote the CCR flow associated to $V$ by $\alpha^{V}$. 

The product system of $\alpha^{V}$, denote it by $F:=\{F(x)\}_{x \in \Omega}$, has the following description. For $x \in \Omega$, $F(x)=\Gamma(Ker(V_x^{*}))$. For a Hilbert space $\clh$, $\Gamma(\clh)$ denotes the symmetric Fock space of $\clh$. 
The product rule on $F$ is given by the following formula:
\[
e(\xi)e(\eta)=e(\xi+V_x\eta)\]
for $\xi \in Ker(V_x^{*})$ and $\eta \in Ker(V_y^{*})$. Here $\{e(\xi)\}$ stands for the exponential vectors of the symmetric Fock space. 

With the foregoing notation, we have the following theorem.

\begin{thm}
\label{characterisation of CCR flows}

The product systems $E$ and $F$ are isomorphic. Moreover the isomorphism $\theta$ is given by the family of unitaries $\{\theta_x:F(x) \to E(x)\}_{x \in \Omega}$ such that for $x \in \Omega$ and $u \in D(x)$,
\[
\theta_x\big(e([u]-[e_x])\big)=\frac{u}{\langle u|e_x \rangle}.
\]
\end{thm}

We omit the proof because the main obstacle in the multivariable case is in the construction of the isometric representation $V$ which is already done. With this in hand, the rest of the  details are exactly the same as the arguments of Arveson. There is one subtle point however which  is the measurability of the map $\theta$ for which we need 
the following lemma. Once the following lemma is established, the measurability of the map $\theta$ can be established as in Section 6.6 of \cite{Arveson}.

\begin{lmma}
Let $E$ be a decomposable product system.  Suppose that $E$ has a unit say $(u_x)_{x \in P}$. Fix $a \in \Omega$ and $e \in D(a)$. Then there exists a left coherent section $(e_x)_{x \in P}$ of decomposable vectors which is measurable such that 
$e_a=e$.
\end{lmma}
\textit{Proof.} It suffices to consider the case when $\langle u_a | e \rangle=1$. 
Fix $x \in P$. By the Archimedean principle, there exists $n \geq 1$ such that $na>x$. Since $e^{n}$ is a decomposable vector in $E(na)$ and $\langle e^{n}|u_{na} \rangle=1$, there exists $e_x \in D(x)$, $f_{x} \in D(na-x)$ with
$\langle e_x|u_x \rangle=1$, $\langle f_x|u_{na-x} \rangle=1$ and $e_xf_x=e^{n}$. Moreover the vectors $e_x$ and $f_x$ are unique and independent of $n$.

Clearly $(e_x)_{x \in P}$ is a left coherent section of decomposable vectors and $e_a=e$. We claim that $\{e_{x}\}_{x \in P}$ is measurable. Fix $n \geq 1$. It suffices to show that $\{e_{x}\}_{x \in P}$ is measurable on $\{0 \leq x \leq a\}$. 
For $0 \leq x \leq a$, let $A_{x}:E(x) \to E(a)$ be defined by $A_x(v)=vu_{na-x}$. Note that $e_x=A_{x}^{*}(e^{n})$. 
Choose a family of measurable sections $\{V_{i}\}$ of the product system $E$ such that for each $x \in P$, $\{V_{i}(x)\}$ is an orthonormal basis for $E(x)$. Write 
\[
e^{n}:=\sum_{i,j}\langle e^{n}|V_{i}(x)V_{j}(na-x) \rangle V_{i}(x)V_{j}(na-x)\]
for $0 \leq x \leq a$. 
Applying $A_{x}^{*}$ to the above expression yields
\[
e_{x}=\sum_{i,j} \langle e^{n}|V_{i}(x)V_{j}(na-x) \rangle \langle V_{j}(na-x) |u_{na-x} \rangle V_{i}(x)\]
for $0 \leq x \leq a$. The measurability of $\{e_{x}\}$ on $\{0 \leq x \leq a\}$ is clear from the above expression. This proves the claim and the proof is complete. \hfill $\Box$

Unlike in the $1$-parameter situtation, it is not true in the multiparameter setting that decomposable product systems necessarily have a unit.   When $d \geq 2$, we construct uncountably many examples
of decomposable product systems which are unitless.  Let $P^{*}$ be the dual of $P$ and $\Omega^{*}$ be the interior of $P^{*}$. Then 
\[
\Omega^{*}=\{\lambda \in \mathbb{R}^{d}: \langle \lambda|x \rangle >0, \forall x \in P \backslash \{0\}\}.\]

Denote the $1$-parameter shift semigroup on $L^{2}(0,\infty)$ by $\{S_{t}\}_{t \geq 0}$. 
Let $\lambda, \mu \in \Omega^{*}$ be given. For $a \in P$, let $\chi(a):=\langle \lambda|a \rangle$. For $b \in P$, set $h_b:=1_{(0,\langle \mu|b \rangle)} \in L^{2}(0,\infty)$. 
For $a \in P$, let $V_{a}:=S_{\langle \mu|a \rangle}$. Then $V:=\{V_{a}\}_{a \in P}$ is an isometric representation of $P$ on $L^{2}(0,\infty)$. 
Let $E$ be the product system associated to the CCR flow $\alpha^{V}$. Set $\overline{E}:=E$. We define an associative multiplication on $\overline{E}$ as follows: for $a, b \in P$, $S \in \overline{E}(a)$ and $T \in \overline{E}(b)$, define
\begin{equation}
\label{product rule}
S.T=S W(\chi(a)h_b)T.
\end{equation}
Here $W(.)$ denotes the Weyl operator. 

Note that for a fixed $a$, $\{W(\chi(a)h_b)\}_{b \in P}$ is a gauge cocycle of the CCR flow $\alpha^{V}$. For a fixed $b \in P$, $\{(W(\chi(a)h_b)\}_{a \in P}$ is a semigroup of unitaries. Using these two facts, it is routine to verify that the product given by Eq. \ref{product rule} is well defined and is associative. The other requirements for $\overline{E}$ to be a product system follows from that of $E$. The main result of \cite{Murugan_Sundar_continuous} ensures that there exists an $E_0$-semigroup, unique up to cocycle conjugacy, such that the associated product system is isomorphic to $\overline{E}$. 

 For $a \in P$, let $\overline{D}(a)$ and $D(a)$ be the decomposable vectors of $\overline{E}(a)$ and $E(a)$ respectively.  A moment's reflection on the definitions reveals that $D(a)=\overline{D}(a)$. 
 Moreover for $a \in P$, $b \in P$, $W(\chi(a)h_b)$ maps $D(b)$ onto $D(b)$. Thus for $a,b \in P$, we have $\overline{D}(a).\overline{D}(b)=\overline{D}(a+b)$. Since $D(a)$ is total in $E(a)$, it follows that $\overline{D}(a)$ is total in $\overline{E}(a)$ for every $a \in P$. Consequently, $\overline{E}$ is decomposable. 
 
 \begin{ppsn}
 \label{unitless}
 Keep the foregoing notation. The following are equivalent. 
 \begin{enumerate}
 \item[(1)] The product system $\overline{E}$ has a unit. 
 \item[(2)] There exists $c>0$ such that $\lambda=c\mu$. 
 \end{enumerate}
  \end{ppsn}
  \textit{Proof.} Suppose that $(1)$ holds. Let $(u_a)_{a \in P}$ be a unit of $\overline{E}$. Then for every $a \in P$, $u_a \in \overline{D}(a)=D(a)$. From the description of the decomposable vectors of CCR flows, it follows that for $a \in P$, there exists $\alpha_a \in \mathbb{C} \backslash \{0\}$ and $\xi_{a} \in Ker(V_a^{*})=L^{2}(0,\langle \mu|a \rangle)$ such that $u_a=\alpha_a T^{(a)}_{e(\xi_a)}$. The equality $u_a.u_b=u_{a+b}$ implies that 
  \[
  \alpha_a \alpha_b T^{(a)}_{e(\xi_a)}W(\chi(a)h_b)T^{(b)}_{e(\xi_b)}=\alpha_{a+b}T^{(a+b)}_{e(\xi_{a+b})}.\]
  This has the consequence that 
  \begin{equation}
  \label{equality of a cocycle}
  \xi_a+V_{a}\xi_b+\chi(a)V_{a}h_b=\xi_{a+b}
    \end{equation}
for $a,b \in P$. 

Fix $b_0 \in \Omega$ such that $\langle \mu|b_0 \rangle=1$. Let $a \in P$ be given. Let $\phi$ be a smooth real valued function such that $supp(\phi) \subset (0,\langle \mu|a \rangle)$. Choose $\delta>0$ such that for $0< t<\delta$, $(0,t)$ is disjoint from $supp(\phi)$. 
Calculate as follows to observe that for $0<t < \delta$,  
\begin{align*}
\int_{-\infty}^{\infty} (\chi(tb_0)h_a(x)+\xi_{a}(x))\phi(x+t)&=\int_{0}^{\infty}(t \chi(b_0)h_a(x)+\xi_a(x))\phi(x+t) \\
&=\langle t\chi(b_0)h_a+\xi_a|V_{tb_0}^{*}\phi \rangle \\
&= \langle V_{tb_0}(\chi(tb_0)h_a+\xi_a)|\phi \rangle \\
&= \langle \xi_{tb_0}+V_{tb_0}(\chi(t_0)h_a+\xi_a)|\phi \rangle ~~( \because supp(\xi_{tb_0}) \cap supp(\phi)=\emptyset)\\
&= \langle \xi_{tb_0+a}|\phi \rangle \\
&=\langle \xi_{a}+V_{a}(\chi(a)h_{tb_0}+\xi_{tb_0})|\phi \rangle \\
&= \langle \xi_{a}|\phi \rangle ~~(\textrm{since $ V_{a}^{*}\phi=0$})\\
&= \int_{-\infty}^{\infty}\xi_{a}(x)\phi(x).
\end{align*}
On rearranging, we get 
\begin{equation}
\label{DE}
\int_{-\infty}^{\infty}\frac{\phi(x+t)-\phi(x)}{t}\xi_a(x)dx= - \langle \lambda|b_0 \rangle \int_{0}^{\langle \mu|a \rangle}\phi(x+t)dx
\end{equation}
Letting $t \to 0$ in Eq. \ref{DE}, we get the differential equation $\xi_{a}^{'}(x)=\langle \lambda|b_0 \rangle$. Here the derivative of $\xi_a$ is in the sense of distributions. Moreover the support of $\xi_a$ is contained in $(0,\langle \mu|a \rangle)$. This implies that there exists $\beta_a \in \mathbb{C}$ such that $\xi_{a}(x)=(\langle \lambda|b_0 \rangle x+\beta_a)1_{(0,\langle \mu|a \rangle)}$. 

Let $a,b \in P$. Now Eq. \ref{equality of a cocycle} implies that for almost all $x \in (0,\langle \mu|a \rangle)$, we have the equality $\langle \lambda|b_0 \rangle x+\beta_a=\langle \lambda|b_0 \rangle x+\beta_{a+b}$. Hence $\beta_{a+b}=\beta_a$ for every $a,b \in P$. Thus $(\beta_a)_{a \in P}$ is constant which we denote by $\beta$. 

Let $a,b \in P$ be given. Eq. \ref{equality of a cocycle} implies that for almost all $x \in (\langle \mu|a \rangle, \langle \mu |a+b \rangle)$, we have the equality
\[
\chi(a)+\langle \lambda|b_0 \rangle (x-\langle \mu|a \rangle)+\beta=\langle \lambda|b_0 \rangle x+\beta.
\]
Thus for $a \in P$, $\langle \lambda|a \rangle= \langle \lambda|b_0 \rangle \langle \mu|a \rangle$. Since $P$ is spanning, it follows that $\lambda=\langle \lambda|b_0 \rangle \mu$. This completes the proof of $(1) \implies (2)$. 

Suppose that $(2)$ holds. Consider the case $P=[0,\infty)$. Denote the $1$-parameter product system $\overline{E}$ defined above for the case when $\mu=1$ and $\lambda=c$ by $\overline{F}$. Now consider the case of a general cone $P$. 
Then $\overline{E}$ is the product system obtained by pulling back the $1$-parameter decomposable product system $\overline{F}$ via the homomorphism $P \ni b \to \langle \mu|b \rangle \in [0,\infty)$. A unit of $\overline{E}$ is given by pulling back a unit of $\overline{F}$. This completes the proof of $(2) \implies (1)$. \hfill $\Box$

It is instructive to work out the isometric representation associated to $\overline{E}$.  For $x \in P$, let  $e_{x}:=T^{(x)}_{e(0)}$. Then $e:=(e_x)_{x \in P}$ is a left coherent section of both $\overline{E}$ and $E$. 
Let $\overline{V}$ be the isometric representation constructed out of the decomposable product system $\overline{E}$. Denote the Hilbert space on which it acts by $H_{\infty}$. We use the notation explained at the beginning of this section. 

It is clear that left coherent sections of $\overline{E}$ and $E$ are the same. This has the consequence that the $e$-logarithm associated to $E$ and $\overline{E}$ are the same. From defintion (see also next section), it follows that for $a \in P$ and $\xi,\eta \in Ker(V_a^{*})$, 
\[
L^{e}(a,T^{(a)}_{e(\xi)},T^{(a)}_{e(\eta)})=\langle \xi|\eta \rangle.\]
The above equation has the consequence that there exists a unitary $U_{a}:ker(V_a^{*}) \to H_a$ such that 
\[
U(\xi-\eta)=[T^{(a)}_{e(\xi)}]-[T^{(a)}_{e(\eta)}]\]
for $a \in \Omega$, $\xi,\eta \in Ker(V_a)^{*}$. It is routine to verify that the unitaries $\{U_a: a \in \Omega\}$ patch together to define a unitary $U:L^{2}(0,\infty) \to H_{\infty}$ such that for $U=U_a$ on $Ker(V_a^{*})$. We claim that $U$ intertwines $V$ and $\overline{V}$. 

Let $a \in \Omega$ be given. For $b \in \Omega$ and $\xi,\eta \in Ker(V_b)^{*}$, calculate as follows to observe that 
\begin{align*}
\overline{V}_{a}U(\xi-\eta)&=\overline{V}_{a}([T^{(b)}_{e(\xi)}]-[T^{(b)}_{e(\eta)}]) \\
&=[T^{(a)}_{e(0)}.T^{(b)}_{e(\xi)}]-[T^{(a)}_{e(0)}.T^{(b)}_{e(\eta)}] \\
&=[T^{(a+b)}_{e(\chi(a)V_ah_b+V_a \xi)}]-[T^{(a+b)}_{e(\chi(a)V_ah_b+V_a \eta)}]\\
&=U(\chi(a)V_ah_b+V_a\xi-\chi(a)V_a h_b-V_a\eta)\\
&=UV_{a}(\xi-\eta).
\end{align*}
Hence $U$ intertwines $V$ and $\overline{V}$. 

In what follows, we denote the product system $\overline{E}$ by $\overline{E}_{\lambda,\mu}$ to stress the dependence of $\overline{E}$ on $\lambda$ and $\mu$. 

\begin{ppsn}
\label{uncountable}
Let $\lambda,\mu_1,\mu_2 \in \Omega^{*}$ be given. Suppose $\overline{E}_{\lambda,\mu_1}$ and $\overline{E}_{\lambda,\mu_2}$ are isomorphic. Then $\mu_1$ and $\mu_2$ are scalar multiplies of each other. 
\end{ppsn}
\textit{Proof.} For $i=1,2$, let $\overline{V}^{(i)}$ be the isometric representation constructed out of the decomposable product system $\overline{E}_{\lambda,\mu_i}$. Then $\overline{V}^{(1)}$ and $\overline{V}^{(2)}$ are unitarily equivalent. 
Fix $i\in \{1,2\}$, the isometric representation $\overline{V}^{(i)}$ is the pull back of the shift semigroup $\{S_t\}_{t \geq 0}$ via the semigroup homomorphism $P \ni a \to \langle \mu_i|a \rangle \in [0,\infty)$. Consequently, the minimal unitary dilation of $\overline{V}^{(i)}$, denote it by $\overline{U}^{(i)}$, is the pull back of the regular representation of $\mathbb{R}$ on $L^{2}(\mathbb{R})$ via the group homomorphism $\mathbb{R}^{d} \ni x \to \langle \mu_i | x \rangle \in \mathbb{R}$. 

But $\overline{U}^{(1)}$ and $\overline{U}^{(2)}$ are unitarily equivalent. Hence 
\[
\{x \in \mathbb{R}^{d}: \langle \mu_1|x \rangle=0 \}= Ker(\overline{U}^{(1)})=Ker(\overline{U}^{(2)})=\{x \in \mathbb{R}^{d}: \langle \mu_2|x \rangle=0\} .
\]
This shows that $\mu_1$ and $\mu_2$ are scalar multiplies of each other. This completes the proof. 

\begin{rmrk}
Thus, in the higher dimensional case, i.e. when $d \geq 2$, Prop. \ref{unitless} coupled with Prop. \ref{uncountable} implies that there are uncountably many decomposable product systems which do not admit a unit.
In the higher dimensional case, when we talk of units, cohomological considerations need to be taken into account (see \cite{Anbu}). The proof of Prop. \ref{unitless} shows that the product system constructed does not admit such``twisted" units either. 

\end{rmrk}

\section{Injectivity of the CCR functor}
In this section, we show that for a pure isometric representation $V$, the cocycle conjugacy class of $\alpha^{V}$ remembers the isometric representation up to unitary equivalence. We prove this by showing that the isometric representation constructed out of the decomposable $E_0$-semigroup $\alpha^{V}$ is $V$. 

Let $V:=(V_{x})_{x \in P}$ be a pure isometric representation on a separable Hilbert space $\clh$. Denote the CCR flow associated to $V$ acting on the algebra of bounded operators on the symmetric Fock space $\Gamma(\clh)$ by $\alpha^{V}=\{\alpha_{x}\}_{x \in P}$. Let $E:=\{E(x)\}_{x \in P}$ be the associated product system. For $a \in \Omega$ and $\xi \in Ker(V_{a}^{*})$, let $T^{(a)}_{e(\xi)}$ be the exponential vectors defined on $\Gamma(\clh)$. Recall that
\[
T^{(a)}_{e(\xi)}e(\eta)=e(\xi+V_{a}\eta)\]
for $\xi \in Ker(V_{a}^{*})$ and $a \in \Omega$. Note that 
\begin{enumerate}
\item[(1)] for $a, b \in \Omega$ and $\xi \in Ker(V_a^{*})$, $T^{(a+b)}_{e(\xi)}=T^{(a)}_{e(\xi)}\Gamma(V_b)$, and
\item[(2)] for $a,b \in \Omega$ and $\xi \in Ker(V_{a}^{*})$, $T^{(a+b)}_{e(V_b\xi)}=\Gamma(V_b)T^{(a)}_{e(\xi)}$.

\end{enumerate}

Let $\widetilde{V}:=(\widetilde{V}_{a})_{a \in \Omega}$ be the semigroup of isometries constructed out of the decomposable $E_0$-semigroup $\alpha^{V}$. Denote the Hilbert space on which $\widetilde{V}$ acts by $H_{\infty}$. 

\begin{ppsn}
\label{Isometric}
With the foregoing notation, we have the following. There exists a unitary $U:\clh \to H_{\infty}$ such that \[UV_{a}U^{*}=\widetilde{V}_{a}\] for every $a \in \Omega$. 
\end{ppsn}
\textit{Proof.} Note that for $a \in \Omega$, the set of decomposable vectors in $E(a)$ is given by \[D(a):=\{\lambda T^{(a)}_{e(\xi)}:\lambda \in \mathbb{C}\backslash \{0\}, \xi \in Ker(V_{a}^{*})\}.\] For $x \in P$, let $e_{x}=\Gamma(V_x)$. Then $e:=(e_x)_{x \in P}$ is a  left coherent section of decomposable vectors of norm $1$.  It follows from definition that \[L^{e}(a,T^{(a)}_{e(\xi)},T^{(a)}_{e(\eta)})=\langle \xi|\eta \rangle\]
for $a \in \Omega$ and $\xi,\eta \in Ker(V_{a}^{*})$.

Calculate as follows to observe that for $a,b \in \Omega$, $\xi \in Ker(V_a^{*})$ and $\eta \in Ker(V_b^{*})$, 
\begin{align*}
\label{equality}
\Big\langle [T^{(a)}_{e(\xi)}]-[T^{(a)}_{e(0)}]\Big|[T^{(b)}_{e(\eta)}]-[T^{(b)}_{e(0)}] \Big\rangle&=\Big\langle [T^{(a)}_{e(\xi)}\Gamma(V_b)]-[T^{(a)}_{e(0)}\Gamma(V_b)]\Big| [T^{(b)}_{e(\eta)}\Gamma(V_a)]-[T^{(b)}_{e(0)}\Gamma(V_a)] \Big \rangle\\
                             &=\Big \langle [T^{(a+b)}_{e(\xi)}]-[T^{(a+b)}_{e(0)}]\Big|[T^{(a+b)}_{e(\eta)}]-[T^{(a+b)}_{e(0)}]\Big\rangle \\
                             &=L^{e}(a+b,T^{(a+b)}_{e(\xi)},T^{(a+b)}_{e(\eta)})-L^{e}(a+b,T^{(a+b)}_{e(\xi)},T^{(a+b)}_{e(0)})\\
                             &~~~-L^{e}(a+b,T^{(a+b)}_{e(0)},T^{(a+b)}_{e(\eta)})+L^{e}(a+b,T^{(a+b)}_{e(0)},T^{(a+b)}_{e(0)}) \\
                             &=\langle \xi|\eta \rangle.
\end{align*}
Note  that $\bigcup_{a \in \Omega}Ker(V_a^{*})$ is total in $\clh$ (for the representation $V$ is pure), and the set $\{[T^{(a)}_{e(\xi)}]-[T^{(a)}_{e(0)}]: a \in \Omega, \xi \in Ker(V_{a}^{*})\}$ is total in $H_{\infty}$. Hence the previous calculation implies that 
there exists a unique unitary $U:\clh \to H_{\infty}$ such that for $a \in \Omega$, $\xi \in Ker(V_{a}^{*})$, 
\[
U\xi=[T^{(a)}_{e(\xi)}]-[T^{(a)}_{e(0)}].\]
 
 We claim that $U$ intertwines $V$ and $\widetilde{V}$. Let $a \in \Omega$ be given. Fix $b \in \Omega$ and $\xi \in Ker(V_{b}^{*})$.  Calculate as follows to observe that
 \begin{align*}
 UV_{a}\xi&=[T^{(a+b)}_{e(V_a\xi)}]-[T^{(a+b)}_{e(0)}] \\
                &=[\Gamma(V_a)T^{(b)}_{e(\xi)}]-[\Gamma(V_a)T^{(b)}_{e(0)}] \\
                &=\widetilde{V}_{a}([T^{(b)}_{e(\xi)}]-[T^{(b)}_{e(0)}])\\
                &=\widetilde{V}_{a}U\xi.
  \end{align*}
  Since $\displaystyle \bigcup_{b \in \Omega}Ker(V_{b}^{*})$ is dense in $\clh$, it follows that $UV_{a}U^{*}=\widetilde{V}_{a}$ for every $a \in \Omega$.  \hfill $\Box$
 
 The following theorem is an immediate consequence of Prop. \ref{Isometric}.
 \begin{thm}
 \label{Injectivity of CCR}
 Let $V$ and $W$ be pure isometric representations of $P$. Denote the CCR flows associated to $V$ and $W$ by $\alpha^{V}$ and $\alpha^{W}$ respectively. Then the following are equivalent.
 \begin{enumerate}
 \item[(1)] The CCR flows $\alpha^{V}$ and $\alpha^{W}$ are cocycle conjugate.
 \item[(2)] The isometric representations $V$ and $W$ are unitarily equivalent. 
  \end{enumerate}
 \end{thm}
 
   Theorem \ref{Injectivity of CCR} gives a more conceptual explanation of the result obtained in \cite{Anbu_Sundar}. Let us quickly recall the main result of  \cite{Anbu_Sundar}. 
  By a $P$-module, we mean a non-empty proper closed subset of $\mathbb{R}^{d}$ say $A$ which is $P$-invariant, i.e. $P+A \subset A$. Fix $k \in \{1,2\cdots,\} \cup \{\infty\}$ and let $\clk$ be a 
  Hilbert space of dimension $k$. For $x \in P$, let $V_{x}$ be the isometry on $L^{2}(A,\clk)$ defined as follows: for $f \in L^{2}(A,\clk)$, 
\begin{equation}
\label{isometries}
V_{x}(f)(y):=\begin{cases}
 f(y-x)  & \mbox{ if
} y -x \in A,\cr
   &\cr
    0 &  \mbox{ if } y-x \notin A.
         \end{cases}
\end{equation}
Then $(V_{x})$ is an isometric representation of $P$ which we denote by $V^{(A,k)}$ and we denote the corresponding CCR flow by $\alpha^{(A,k)}$. 

The main result obtained in  \cite{Anbu_Sundar} is the
following. Let $A_1,A_2$ be two $P$-modules and $k_1,k_2 \in \cdots\{1,2,\cdots,\} \cup \{\infty\}$ be given. Then the following are equivalent.
\begin{enumerate}
\item[(1)] The CCR flows $\alpha^{(A_1,k_1)}$ and $\alpha^{(A_2,k_2)}$ are  cocycle conjugate.
\item[(2)] The modules $A_1$ and $A_2$ are translates of each other and $k_1=k_2$.
\end{enumerate}
 In \cite{Anbu_Sundar}, the equivalence between $(1)$ and $(2)$ was established using groupoid machinery after a  careful analysis of associated gauge groups and after overcoming two problems regarding the
unitary groups of von Neumann algebras. (See Section 4 of \cite{Anbu_Sundar}). Now there is a more efficient way of achieving the equivalence between $(1)$ and $(2)$. Let us consider a third statement
\begin{enumerate}
\item[(3)] The representations $V^{(A_1,k_1)}$ and $V^{(A_2,k_2)}$ are unitarily equivalent. 
\end{enumerate}  
Theorem \ref{Injectivity of CCR} asserts that $(1)$ and $(3)$ are equivalent. It is clear that $(2)$ implies $(3)$. The point we wish to stress is that it is only in the proof of $(3) \implies (2)$ that we need groupoids.

 Keep the foregoing notation. Suppose $(3)$ holds.
Assume that $A_1$ is not a translate of $A_2$. The calculation carried out in Prop. 4.7 of \cite{Anbu_Sundar} asserts that $V^{(A_1,k_1)}$ and $V^{(A_2,k_2)}$ are disjoint which is a contradiction. Consequently, $A_1$ and $A_2$ are translates of each other.
By Corollary 3.4 of \cite{Anbu_Sundar}, we have that for $i=1,2$, the commutant of the von Neumann algebra generated by $\{V^{(A_i,k_i)}_{x}:x \in P\}$ is of the form $M_i \otimes B(\clk_i)$ where
$M_i$ is an abelian von Neumann algebra. Choose an isomorphism say $\phi:M_1 \otimes B(\clk_1) \to M_2 \otimes B(\clk_2)$. Suppose $k_1>k_2$.  Consider a character $\chi$ of the commutative $C^{*}$-algebra $M_2$. Then the map
\[
B(\clk_1) \ni T \to (\chi \otimes 1)(\phi(T)) \in B(\clk_2)\]
is a non-zero $*$-homomorphism which is a contradiction, since $B(\clk_1)$ does not admit a non-zero representation on a Hilbert space whose dimension is strictly less than $\dim(\clk_1)$. Hence $k_1\leq k_2$. By symmetry, it follows that $k_2 \leq k_1$ and thus $k_1=k_2$. This 
completes the proof. 

\section{Local cocycles for CCR flows}
 In the one parameter case, unitary local cocycles and positive contractive local cocyles for CCR flows were computed by Arveson (\cite{Arv_Fock}) and Bhat (\cite{Bhat_cocycles}) respectively. The analysis once again relies on the fact that there are enough units for one parameter CCR flows.  The situation in the higher dimensional case is different and the examples that we know so far admits only one unit. However with the explicit knowledge of decomposable vectors and  the $e$-logarithm in hand, we could compute the local cocycles of multiparameter CCR flows. In this section, as an application of our construction, we derive a formula for  the positive contractive local cocyles of a CCR flow associated to a pure isometric representation. The local unitary cocyles were already computed in \cite{Anbu}.

Let us recall the definition of a local cocycle. 
\begin{dfn}
\label{local cocycle}
Let $\alpha:=\{\alpha_{x}\}_{x \in P}$ be an $E_0$-semigroup on $B(\clh)$. A family $C:=\{C_{x}\}_{x \in P}$ of bounded operators on $\clh$ is said to be a local cocycle for $\alpha$ if 
\begin{enumerate}
\item[(1)] the map $P \ni x \to C_{x} \in B(\clh)$ is weakly continuous and $C_{0}=Id$,
\item[(2)] for $x,y \in P$, we have the cocycle condition $C_{x+y}=C_{x}\alpha_{x}(C_y)$, and
\item[(3)] for $x \in P$, $C_{x} \in \alpha_{x}(B(\clh))^{'}$.
\end{enumerate}
The local cocycle $\{C_{x}\}_{x \in P}$ is said to be a unitary (contractive, positive, projective) local cocycle if each $C_{x}$ is unitary (contractive, positive, projective). 
\end{dfn}
\begin{rmrk}
A theorem of Powers which asserts that $E$-semigroups subordinate to a given $E_0$-semigroup $\alpha$ are in one-one correspondence with the set of projective local cocycles of $\alpha$ stays true in the
multiparameter context as well. Similarly, positive contractive local cocycles are in one-one correspondence with quantum dynamical semigroups which are dominated by the given $E_0$-semigroup. For more details, we 
refer the reader to \cite{Powers} and \cite{Bhat_cocycles}.
\end{rmrk}

Let $\alpha:=\{\alpha_{x}\}_{x \in P}$ be an $E_0$-semigroup and let $C:=\{C_{x}\}_{x \in P}$ be a family of bounded operators on $\clh$ such that $C_{x} \in \alpha_{x}(B(\clh))^{'}$. Denote the product system associated to $\alpha$ by $E:=\{E(x)\}_{x \in P}$. 
For $x \in P$, define $\theta_{x}:E(x) \to E(x)$ by 
\begin{equation}
\label{morphism}
\theta_{x}(T)=C_{x}T
\end{equation}
for $T \in E(x)$. Note that $\theta_{x}$ is well defined and $||\theta_x||=||C_x||$. It is routine to verify that the cocycle condition $C_x\alpha_x(C_y)=C_{x+y}$ is equivalent to the following condition: 
for $x,y \in P$, $T \in E(x)$ and $S \in E(y)$, $\theta_{x+y}(TS)=\theta_x(T)\theta_y(S)$. In other words, $\{C_x\}_{x \in P}$ satisfies the cocycle condition if and only $\theta:=\{\theta_x\}_{x \in P}$ is a morphism of the product system $E$. One consequence of seeing local cocycles this way is that the adjoint of a local cocycle is again a local cocycle.

\begin{ppsn}
\label{measurability}
Let $\alpha:=\{\alpha_x\}_{x \in P}$ be an $E_0$-semigroup on $B(\clh)$ and $C:=\{C_x\}_{x\in P}$ be a family of bounded operators on $\clh$ satisfying Conditions $(2)$ and $(3)$ with $C_0=Id$.   Assume that for some $a \in \Omega$, $C_{a} \neq 0$. Suppose in addition that
we have the following measurability hypothesis i.e. 
\begin{enumerate}
\item[$(1)^{'}$] for $\xi,\eta \in \clh$, the map $P \ni x \to \langle C_x\xi|\eta \rangle \in \mathbb{C}$ is measurable.
\end{enumerate}
Then $C:=\{C_x\}_{x \in P}$ is a local cocycle i.e. it satisfies Condition $(1)$ of Defn. \ref{local cocycle}. 
\end{ppsn}
\textit{Proof.} Let $\theta:=\{\theta_{x}\}_{x \in P}$ be the morphism, given by Eq. \ref{morphism}, corresponding to $\{C_x\}_{x \in P}$ of the associated product system $E$. 
Observe that if we identify the Hilbert space $E(x+y)$ with the tensor product $E(x) \otimes E(y)$, then $\theta_{x+y}=\theta_{x} \otimes \theta_{y}$. Hence $||\theta_{x+y}|| = ||\theta_x||||\theta_y||$ for $x,y \in P$.
 For $x \in P$, let $f:P \to [0,\infty)$ be defined by
\[f(x)=||\theta_x||=||C_x||\] for $x \in P$. Then $f(x+y)=f(x)f(y)$ for $x,y \in P$ and $f(na) \neq 0$ for every $n \geq 1$. By the Archimedean property, it follows that $f(x) \neq 0$ for every $x \in P$. Also $f$ is measurable. Hence there exists $\lambda \in \mathbb{R}^{d}$ such that 
$f(x)=e^{\langle \lambda|x\rangle}$.

For $x \in P$, let $\widetilde{C}_{x}=e^{-\langle \lambda|x \rangle}C_{x}$. Note that $\{\widetilde{C}_{x}\}_{x \in P}$ satisfies the hypothesis of the proposition. Moreover $||\widetilde{C}_x||=1$ for $x \in P$. It suffices to prove that $\{\widetilde{C}_{x}\}_{x \in P}$ is weakly continuous. 
To see this, define for $x \in P$, a linear map $\rho_{x}$ on $B(\clh)$ by the formula
\[
\rho_{x}(T)=\widetilde{C}_{x}\alpha_{x}(T)
\]
for $T \in B(\clh)$. Then $\{\rho_{x}\}_{x \in P}$ forms a semigroup of linear maps on $B(\clh)$ which is weakly measurable. Moreover $\{\rho_{x}\}_{x \in P}$ is uniformly bounded. 

We claim that $\displaystyle \bigcap_{x \in \Omega}ker(\rho_{x})=\{0\}$. Suppose that there exists $A \in B(\clh)$ such that $\rho_{x}(A)=0$ for all $x \in \Omega$. Fix $x \in \Omega$. Then $\widetilde{C}_{x}\alpha_{x}(A)=0$. But $\widetilde{C}_{x}$ commutes with $\alpha_{x}(B(\clh))$. Hence \[\widetilde{C}_{x}\alpha_{x}(XAY)=0\] for every $X,Y \in B(\clh)$. Since $B(\clh)$ does not have any nontrivial $\sigma$-weakly closed ideal and $\widetilde{C}_x \neq 0$, it follows that $A=0$. This proves our claim.

By adapting the  proof of 
Prop. 4.2 of \cite{Murugan_Sundar}, we see that for $T \in B(\clh)$ and $\xi,\eta \in \clh$, the map $P \ni x \to \langle \rho_x(T)\xi|\eta \rangle \in \mathbb{C}$ is continuous. Taking $T=Id$, we obtain that 
for $ \xi,\eta \in \clh$, the map $P \ni x \to \langle \widetilde{C}_x\xi|\eta \rangle \in \mathbb{C}$ is continuous. This completes the proof. \hfill $\Box$

Fix a decomposable $E_0$-semigroup $\alpha:=\{\alpha_x\}_{x \in P}$ and let $E:=\{E(x)\}_{x \in P}$ be the associated product system.  
First note that if $\{u_x\}_{x \in P}$ is a unit of $\alpha$ and $\{C_{x}\}_{x \in P}$ is a local cocycle for $\alpha$, then $\{C_xu_x\}_{x \in P}$ is a unit. 
We need to know that local cocycles of $\alpha$ map decomposable vectors of $E$ to decomposable vectors of $E$. This follows from the next result.
Let $\clk$ be a Hilbert space of dimension $k$ where $k \in \{1,2,\cdots,\} \cup \{\infty\}$ and $\{S_{t}\}_{t \geq 0}$ be the shift semigroup on $L^{2}((0,\infty), \clk)$.
Denote the CCR flow associated to $\{S_{t}\}_{t \geq 0}$ by $\alpha:=\{\alpha_t\}_{t \geq 0}$. With the foregoing notation, we have the following.

\begin{ppsn}
\label{non-vanishing}
Let $C:=\{C_t\}_{t \geq 0}$ be a local cocycle for $\alpha$. Denote the morphism of the associated product system corresponding to $\{C_{t}\}_{t \geq 0}$ by $\theta:=\{\theta_t\}_{t \geq 0}$. Then for every 
$t \geq 0$ and $\xi \in Ker(S_t^{*})$, 
$\theta_{t}T^{(t)}_{e(\xi)} \neq 0$ where $T^{(t)}_{e(\xi)}$ are the exponential vectors defined by Eq. \ref{exponential}.
 \end{ppsn}
 \textit{Proof.} The notation that we use are inspired by Bhat's notation of Theorem 7.5 of \cite{Bhat_cocycles}. 
 Note that for every $y \in \clk$, $\{\theta_{t}(T^{(t)}_{e(1_{(0,t)}\otimes y))}\}_{t \geq 0}$ is a unit. Hence there exists $B_{y} \in \clk$ and $p_{y} \in \mathbb{C}$ such that 
 \[
 \theta_{t}T^{(t)}_{e(1_{(0,t)}\otimes y)}=e^{tp_{y}}T^{(t)}_{e(1_{(0,t)}\otimes B_{y})}.\]
 Similary for every $y \in \clk$, there exists $q_{y} \in \mathbb{C}$ and $D_{y} \in \clk$ such that  \[
 \theta_{t}^{*}T^{(t)}_{e(1_{(0,t)}\otimes y)}=e^{tq_{y}}T^{(t)}_{e(1_{(0,t)}\otimes D_{y})}.\]

Fix $t>0$ and $\xi \in Ker(S_{t}^{*})$ be given. Suppose that $\theta_{t}T^{(t)}_{e(\xi)} = 0$. Let $\ell \geq 1$ be given. Since step functions are dense in $L^{2}((0,t),\clk)$, for every $\ell \geq 1$, there exists $n_{\ell} \geq 1$, $y^{(\ell)}_1,y^{(\ell)}_2,\cdots, y^{(\ell)}_{n_\ell} \in \clk$ and a partition $0=s^{(\ell)}_{0}<s^{(\ell)}_{1}<s^{(\ell)}_{2}< \cdots <s^{(\ell)}_{n_\ell}=t$ such that 
\[
\big|\big|\xi-\sum_{i=1}^{n_\ell}1_{(s^{(\ell)}_{i-1},s^{(\ell)}_{i})}\otimes y_{i}^{(\ell)} \big|\big| \leq \frac{1}{\ell}.
\]
Set $\xi_{\ell}=\sum_{i=1}^{n_\ell}1_{(s^{(\ell)}_{i-1},s^{(\ell)}_{i})}\otimes y^{(\ell)}_{i}$.
For $\ell \geq 1$ and $1 \leq i \leq n_\ell$, set $t^{(\ell)}_i=s^{(\ell)}_i-s^{(\ell)}_{i-1}$. Observe that $\sum_{i=1}^{n_\ell}t^{(\ell)}_{i}=t$. Note that as $\ell \to \infty$, \[
\prod_{i=1}^{n_\ell}T^{(t^{(\ell)}_i)}_{e\big(1_{(0,t^{(\ell)}_{i})}\otimes y^{(\ell)}_{i}\big)} \to T^{(t)}_{e(\xi)}.\]
Hence as $\ell \to \infty$, 
\[
\Big \langle\theta_{t}\Big(\prod_{i=1}^{n_\ell}T^{(t^{(\ell)}_i)}_{e\big(1_{(0,t^{(\ell)}_{i})}\otimes y^{(\ell)}_{i}\big)}\Big)\Big|e(0)\Big\rangle \to 0.\]
Calculate as follows to observe that
\begin{align*}
\Big \langle\theta_{t}\Big(\prod_{i=1}^{n_\ell}T^{(t^{(\ell)}_i)}_{e\big(1_{(0,t^{(\ell)}_{i})}\otimes y^{(\ell)}_{i}\big)}\Big)\Big|e(0)\Big\rangle =& \Big \langle \Big(\prod_{i=1}^{n_\ell}T^{(t^{(\ell)}_i)}_{e\big(1_{(0,t^{(\ell)}_{i})}\otimes y^{(\ell)}_{i}\big)}\Big)\Big|\theta_{t}^{*}e(0)\Big\rangle \\
&= \Big \langle \Big(\prod_{i=1}^{n_\ell}T^{(t^{(\ell)}_i)}_{e\big(1_{(0,t^{(\ell)}_{i})}\otimes y^{(\ell)}_{i}\big)}\Big)\Big|e^{tq_0}\Big(\prod_{i=1}^{n_\ell}T^{(t^{(\ell)}_i)}_{e\big(1_{(0,t^{(\ell)}_{i})}\otimes D_0\big)}\Big)\Big\rangle\\
&=e^{t\overline{q_0}}e^{\sum_{i=1}^{n_\ell}t^{(\ell)}_{i}\langle y^{(\ell)}_{i}|D_0 \rangle}\\
&=e^{t \overline{q_0}}e^{\int_{0}^{t}\langle \xi_{\ell}(s)|D_0 \rangle} \\
& \to e^{t\overline{q_0}}e^{\int_{0}^{t}\langle \xi(s)|D_0\rangle} \neq 0
\end{align*}
which is a contradiction. Hence $\theta_{t}(T^{(t)}_{e(\xi)}) \neq 0$. This completes the proof. \hfill $\Box$

Prop. \ref{non-vanishing} coupled with the fact that $1$-parameter decomposable product systems are isomorphic to CCR flows yields immediately the following. Let $\alpha:=\{\alpha_{x}\}_{x \in P}$ be a decomposable $E_0$-semigroup with product system $E:=\{E(x)\}_{x \in P}$. Consider a local cocycle $C:=\{C_{x}\}_{x \in P}$ of $\alpha$. Denote the morphism of $E$ associated to $C$ by $\theta:=\{\theta_x\}_{x \in P}$.
Let $x \in P$ be given. If $u \in E(x)$ is  decomposable  then $\theta_x(u)$ is also decomposable. Thus $\theta$ maps decomposable vectors to decomposable vectors.

Let $V:=(V_a)_{a \in \Omega}$ be the isometric representation constructed in Section 4 corresponding to the decomposable $E_0$-semigroup $\alpha$. We keep the notation used in the paragraphs preceeding Prop. \ref{separability}. Let $H_{\infty}$ be the Hilbert space on which $V$ acts. Our next proposition states that every local cocycle of $\alpha$ induces a bounded linear operator on $H_{\infty}$ which lies in the commutant of the von Neumann algebra generated by $\{V_{a}:a \in \Omega\}$. This is the conceptual reason behind the appearance of the commutant in the formula of the gauge group of a CCR flow (see Theorem 7.2 of \cite{Anbu}).

Let $C:=\{C_x\}_{x \in P}$ be a local cocycle of $\alpha$ and $\theta:=\{\theta_x\}_{x \in P}$ be the associated morphism of $E$. Let $\lambda \in \mathbb{R}^{d}$ be such that for $x \in P$, $||\theta_x||=e^{\langle \lambda|x\rangle}$.
\begin{ppsn}
\label{appearance of commutant}
Keep the foregoing notation.  Then there exists a unique contraction denoted $\widetilde{\theta}$ on $H_{\infty}$ such that
\begin{enumerate}
\item[(1)] for $a \in \Omega$, $V_{a}\widetilde{\theta}=\widetilde{\theta}V_{a}$ and $V_{a}^{*}\widetilde{\theta}=\widetilde{\theta}V_{a}^{*}$, and
\item[(2)] for $a \in \Omega$, $u,v \in D(a)$, $\widetilde{\theta}([u]-[v])=[\theta_{a}(u)]-[\theta_{a}(v)]$.
\end{enumerate}
\end{ppsn}
Before we take up the proof of Prop. \ref{appearance of commutant}, let us recall the following remarkable formula due to Arveson which expresses the $e$-logarithm in terms of partitions. Fix a left coherent section $e:=\{e_{x}\}_{x \in P}$ of decomposable vectors of unit norm. 

Fix $a \in \Omega$ and let $u,v \in D(a)$ be given. Let \[\mathcal{P}:=\{0=t_0<t_1<t_2<\cdots<t_n=1\}\] be a partition of the unit interval $[0,1]$. Choose $\{u_{i}\}_{i=1}^{n}$ and $\{v_i\}_{i=1}^{n}$ such that for every $i$,  $u_i,v_i \in D((t_i-t_{i-1})a)$, $u=u_1u_2\cdots,u_n$ and $v=v_1v_2\cdots v_n$. Set $e_{i}=e(t_{i-1}a,t_ia)$. Note that $\displaystyle \sum_{i=1}^{n} \Big(\frac{\langle u_i|v_i \rangle}{\langle u_i|e_i \rangle \langle e_i|v_i \rangle}-1\Big)$ is independent of the choice of $\{u_{i}\}$ and $\{v_i\}$. Arveson's remarkable formula is that
\begin{equation}
\label{remarkable}
L^{e}(a,u,v)=\lim_{\mathcal{P}}\sum_{i=1}^{n}\Big(\frac{\langle u_i|v_i \rangle}{\langle u_i|e_i \rangle \langle e_i|v_i \rangle}-1\Big)
\end{equation}
where the limit is taken over all partitions of $[0,1]$. Here the partitions are partially ordered by the usual notion of refinement. (See Section 6.5 of \cite{Arveson} and Remark \ref{comparision}).

\begin{lmma}
\label{boundedness of induced map}
Let $u^{(1)},u^{(2)},\cdots,u^{(n)} \in D(a)$ be given and $\lambda_1,\lambda_2,\cdots,\lambda_n \in \mathbb{C}$ be such that $\sum_{i=1}^{n}\lambda_{i}=0$. 
Then \[
\sum_{i,j=1}^{n}\lambda_{i}\overline{\lambda_j}L^{e}(a,\theta_{a}(u^{(i)}),\theta_{a}(u^{(j)})) \leq \sum_{i,j=1}^{n} \lambda_{i}\overline{\lambda_j}L^{e}(a,u^{(i)},u^{(j)}).
\]

\end{lmma}
\textit{Proof.} Since the $e$-logarithm, $(L^{e}(a,.,.))_{a \in \Omega}$ is homogeneous, by replacing $\{\theta_{x}\}_{x \in P}$ by $\{e^{-\langle \lambda|x \rangle}\theta_{x}\}_{x \in P}$, we can assume that $||\theta_{x}||=1$ for every $x \in P$. 
By restricting attention to the one parameter product system $\{E(ta):t>0\}$, we see that it is enough to consider the case when the dimension of the cone $P$ is $1$. Thus we can take $P=[0,\infty)$ and $a=1$. But $1$-parameter decomposable
$E_0$-semigroups are CCR flows and have units in abundance. In view of Theorem 6.4.2 of \cite{Arveson}, we can take our left coherent section $(e_{t})_{t \geq 0}$ to be a unit. 

Set $f_{t}:=\theta_{t}^{*}e_{t}$. Then $(f_{t})_{t \geq 0}$ is a unit. Choose $c \in \mathbb{R}$ such that $f_{t}^{*}f_{t}=e^{ct}$. Define $\widetilde{f_t}=\frac{f_t}{||f_t||}$. For each $i$, choose a left coherent section $(u^{(i)})_{t \geq 0}$ such that $u^{(i)}_{1}=u^{(i)}$. Let $\epsilon>0$ be given. Choose $\delta>0$ such that if $|t|<\delta$ then $e^{-ct}<1+\epsilon$.

Let $\mathcal{P}:=\{0=t_0<t_1<t_2<\cdots<t_m=1\}$ be a partition of $[0,1]$. Denote the norm or the mesh of $\mathcal{P}$ by $||\mathcal{P}||$. For $i \in \{1,2,\cdots,m\}$, let $s_i=t_i-t_{i-1}$.
We claim that if $||\mathcal{P}|| < \delta$ then
\begin{align}
\label{inequality}
&\sum_{i,j=1}^{n}\lambda_{i}\overline{\lambda_j}\Big(\sum_{k=1}^{m}\big(\frac{\langle \theta_{s_k}u^{(i)}(t_{k-1},t_k)|\theta_{s_k}u^{(j)}(t_{k-1},t_k)\rangle}{\langle \theta_{s_k}u^{(i)}(t_{k-1},t_k)|e(t_{k-1},t_k)\rangle \langle e(t_{k-1},t_k)|\theta_{s_k}u^{(j)}(t_{k-1},t_k)\rangle}-1\big) \Big) & \\
&\leq (1+\epsilon)\sum_{i,j=1}^{n}\lambda_{i}\overline{\lambda_j}\Big(\sum_{k=1}^{m}\big(\frac{\langle u^{(i)}(t_{k-1},t_k)|u^{(j)}(t_{k-1},t_k)\rangle}{\langle u^{(i)}(t_{k-1},t_k)|\widetilde{f}(t_{k-1},t_k)\rangle \langle \widetilde{f}(t_{k-1},t_k)|u^{(j)}(t_{k-1},t_k)\rangle}-1\big) \Big).
\end{align}
Fix $k \in \{1,2,\cdots,m\}$. Note that the matrix $(\langle \theta_{s_k}u^{(i)}(t_{k-1},t_k)|\theta_{s_k}u^{(j)}(t_{k-1},t_k)\rangle)$ is positive and since $\theta_{s_k}$ is a contraction, we have \[(\langle \theta_{s_k}u^{(i)}(t_{k-1},t_k)|\theta_{s_k}u^{(j)}(t_{k-1},t_k)\rangle) \leq (\langle u^{(i)}(t_{k-1},t_k)|u^{(j)}(t_{k-1},t_k)\rangle).\] Thus given $\mu_1,\mu_2,\cdots,\mu_n \in \mathbb{C}$, 
\begin{equation}
\label{positivity}
\sum_{i,j}\mu_i\overline{\mu_j}\langle \theta_{s_k}u^{(i)}(t_{k-1},t_k)|\theta_{s_k}u^{(j)}(t_{k-1},t_k)\rangle \leq \sum_{i,j}\mu_i\overline{\mu_j}\langle u^{(i)}(t_{k-1},t_k)|u^{(j)}(t_{k-1},t_k)\rangle .
\end{equation}
Calculate as follows to observe that 
\begin{align*}
&\sum_{i,j=1}^{n}\lambda_{i}\overline{\lambda_j}\Big(\sum_{k=1}^{m}\big(\frac{\langle \theta_{s_k}u^{(i)}(t_{k-1},t_k)|\theta_{s_k}u^{(j)}(t_{k-1},t_k)\rangle}{\langle \theta_{s_k}u^{(i)}(t_{k-1},t_k)|e_{t_k-t_{k-1}}\rangle \langle e_{t_k-t_{k-1}}|\theta_{s_k}u^{(j)}(t_{k-1},t_k)\rangle}-1\big) \Big)  \\
&=\sum_{k=1}^{m} \Big(\sum_{i,j=1}^{n}\lambda_{i}\overline{\lambda_j}\big(\frac{\langle \theta_{s_k}u^{(i)}(t_{k-1},t_k)|\theta_{s_k}u^{(j)}(t_{k-1},t_k)\rangle}{\langle \theta_{s_k}u^{(i)}(t_{k-1},t_k)|e_{t_k-t_{k-1}}\rangle \langle e_{t_k-t_{k-1}}|\theta_{s_k}u^{(j)}(t_{k-1},t_k)\rangle}-1 \big)\Big)\\
&=\sum_{k=1}^{m} \Big(\sum_{i,j=1}^{n}\lambda_{i}\overline{\lambda_j}\frac{\langle \theta_{s_k}u^{(i)}(t_{k-1},t_k)|\theta_{s_k}u^{(j)}(t_{k-1},t_k)\rangle}{\langle \theta_{s_k}u^{(i)}(t_{k-1},t_k)|e_{t_k-t_{k-1}}\rangle \langle e_{t_k-t_{k-1}}|\theta_{s_k}u^{(j)}(t_{k-1},t_k)\rangle}\Big)\\
\end{align*}
\begin{align*}
&\leq \sum_{k=1}^{m} \Big(\sum_{i,j=1}^{n}\lambda_{i}\overline{\lambda_j}\frac{\langle u^{(i)}(t_{k-1},t_k)|u^{(j)}(t_{k-1},t_k)\rangle}{\langle \theta_{s_k}u^{(i)}(t_{k-1},t_k)|e_{t_k-t_{k-1}}\rangle \langle e_{t_k-t_{k-1}}|\theta_{s_k}u^{(j)}(t_{k-1},t_k)\rangle}\Big) ~(\text{by Eq. \ref{positivity}}) \\
& \leq \sum_{k=1}^{m} \Big(\sum_{i,j=1}^{n}\lambda_{i}\overline{\lambda_j}\big(\frac{\langle u^{(i)}(t_{k-1},t_k)|u^{(j)}(t_{k-1},t_k)\rangle}{\langle u^{(i)}(t_{k-1},t_k)|||f_{t_k-t_{k-1}}||\widetilde{f}_{t_k-t_{k-1}}\rangle \langle ||f_{t_k-t_{k-1}}||\widetilde{f}_{t_k-t_{k-1}}|u^{(j)}(t_{k-1},t_k)\rangle}-1\big)\Big)\\
& \leq \sum_{k=1}^{m}e^{c(t_{k-1}-t_k)} \Big(\sum_{i,j=1}^{n}\lambda_{i}\overline{\lambda_j}\big(\frac{\langle u^{(i)}(t_{k-1},t_k)|u^{(j)}(t_{k-1},t_k)\rangle}{\langle u^{(i)}(t_{k-1},t_k)|\widetilde{f}_{t_k-t_{k-1}}\rangle \langle \widetilde{f}_{t_k-t_{k-1}}|u^{(j)}(t_{k-1},t_k)\rangle}-1\big)\Big) \\
& \leq (1+\epsilon)\sum_{i,j=1}^{n}\lambda_{i}\overline{\lambda_j}\Big(\sum_{k=1}^{m}\big(\frac{\langle u^{(i)}(t_{k-1},t_k)|u^{(j)}(t_{k-1},t_k)\rangle}{\langle u^{(i)}(t_{k-1},t_k)|\widetilde{f}_{t_k-t_{k-1}}\rangle \langle \widetilde{f}_{t_k-t_{k-1}}|u^{(j)}(t_{k-1},t_k)\rangle}-1\big) \Big).
\end{align*}
This proves our claim. Taking limit in Eq. \ref{inequality} over the directed set of partitions with mesh less than $\delta$, we obtain
\[
\sum_{i,j=1}^{n}\lambda_{i}\overline{\lambda_j}L^{e}(1,\theta_{1}(u^{(i)}),\theta_{1}(u^{(j)})) \leq (1+\epsilon)\sum_{i,j=1}^{n} \lambda_{i}\overline{\lambda_j}L^{\widetilde{f}}(1,u^{(i)},u^{(j)}).
\]
Letting $\epsilon \to 0$, we have
\[
\sum_{i,j=1}^{n}\lambda_{i}\overline{\lambda_j}L^{e}(1,\theta_{1}(u^{(i)}),\theta_{1}(u^{(j)})) \leq \sum_{i,j=1}^{n} \lambda_{i}\overline{\lambda_j}L^{\widetilde{f}}(1,u^{(i)},u^{(j)}).
\]

 But Theorem 6.4.2 of \cite{Arveson} implies that \[\sum_{i,j=1}^{n} \lambda_{i}\overline{\lambda_j}L^{\widetilde{f}}(1,u^{(i)},u^{(j)})=\sum_{i,j=1}^{n} \lambda_{i}\overline{\lambda_j}L^{e}(1,u^{(i)},u^{(j)}).\]  The proof is now complete. \hfill $\Box$
 
\textit{Proof of Prop. \ref{appearance of commutant}.} Fix $a \in \Omega$. Consider the vector space $\mathbb{C}_{0}\Delta(a)$ together with its semi-definite inner product. The estimate obtained in Lemma \ref{boundedness of induced map} implies that the map
\[\mathbb{C}_{0}\Delta(a) \ni \sum_{u}f(u)\delta_u \to \sum_{u}f(u)\delta_{\theta_{a}(u)} \in \mathbb{C}_0\Delta(a)\]
descends to a contraction on $H_{a}$. It is also clear that the maps so obtained on $H_{a}$ patch together to define a well defined contraction on $H_{\infty}$ which we denote by $\widetilde{\theta}$. It is clear that
\[
\widetilde{\theta}([u]-[v])=[\theta_a(u)]-[\theta_{a}(v)]
\]
for $a \in \Omega$ and $u,v \in D(a)$. Note that $\widetilde{\theta}$ leaves $Ker(V_{a}^{*})=\{[u]-[v]: u,v \in D(a)\}$ invariant for every $a \in \Omega$. Let $a, b \in \Omega$. Choose $e \in D(a)$. Calculate as follows to observe that for $u,v \in D(b)$
\begin{align*}
V_{a}\widetilde{\theta}([u]-[v])&=V_{a}([\theta_b(u)]-[\theta_{b}(v)]) \\
                                               &=[\theta_{a}(e)\theta_{b}(u)]-[\theta_a(e)\theta_{b}(v)] \\
                                               &=[\theta_{a+b}(eu)]-[\theta_{a+b}(ev)] \\
                                               &=\widetilde{\theta}([eu]-[ev]) \\
                                               &=\widetilde{\theta}V_{a}([u]-[v]).
\end{align*}
This proves that $\widetilde{\theta}$ commutes with $\{V_{a}: a \in \Omega\}$. 

Fix $a \in \Omega$. Set $S:=V_{a}$. By the Archimedean property and the purity of the representation $\{V_{b}\}_{b \in \Omega}$, it follows that  $S$ is a pure isometry. Moreover $\widetilde{\theta}$ commutes with $\{S^{n}: n \in \mathbb{N}\}$ and leaves $Ker(S^{*})$ invariant. Hence $\widetilde{\theta}$ commutes with $S^{*}$. (This could be argued for instance by considering $S$ to be the standard shift on $\ell^{2}(\mathbb{N})$ with multiplicity which is possible due to Wold decomposition. We leave the details to the reader.)
Consequently, $V_{a}^{*}\widetilde{\theta}=\widetilde{\theta}V_{a}^{*}$. The uniqueness part is obvious. This completes the proof. \hfill $\Box$

Keep the foregoing notation. Let $\theta^{*}=\{\theta_{x}^{*}\}_{x \in P}$ be the morphism of $E$ associated to the local cocycle $\{C_{x}^{*}\}_{x \in P}$. Fix a left coherent section $e:=(e_x)_{x \in P}$ of decomposable vectors with 
unit norm. 

\begin{lmma}
\label{adjoint formula}
Consider an element $a \in \Omega$. Let $\lambda_1,\lambda_2,\cdots,\lambda_{m}, \mu_1,\mu_2,\cdots,\mu_n \in \mathbb{C}$ be such that 
\[
\sum_{i=1}^{m}\lambda_i=\sum_{j=1}^{n}\mu_j=0.\]
Suppose $u^{(1)},u^{(2)},\cdots,u^{(m)}, v^{(1)},v^{(2)},\cdots,v^{(n)} \in D(a)$. Then
\[
\sum_{i,j}\lambda_{i}\overline{\mu_j}L^{e}(a,\theta_{a}u^{(i)},v^{(j)})=\sum_{i,j}\lambda_{i}\overline{\mu_j}L^{e}(a,u^{(i)},\theta_{a}^{*}v^{(j)}).\]
\end{lmma}
\textit{Proof.} Taking Remark \ref{comparision} into account, we see that the assertion is really  one parameter in nature. Thus, we can assume $P=[0,\infty)$. Theorem 6.4.2 of \cite{Arveson} ensures that the sum $\sum_{i,j}\lambda_i \overline{\mu_j}L^{e}(t,\theta_{t}u^{(i)},v^{(j)})$ is independent of the chosen section $e$. Consequently we can choose $e$ to be a unit. Set $\widetilde{f}_{t}=\theta_{t}e_{t}$ and $\widetilde{g}_{t}=\theta_{t}^{*}e_{t}$. Let $f_{t}=\frac{\widetilde{f}_{t}}{||\widetilde{f}_{t}||}$ and $g_{t}=\frac{\widetilde{g}_{t}}{||\widetilde{g}_{t}||}$. Note that $(\widetilde{f}_{t}), (\widetilde{g}_{t}), (f_{t}),(g_{t})$ are units.
Let $\alpha,\beta,\gamma,\delta \in \mathbb{C}$ be such that $||\widetilde{f}_{t}||=e^{\alpha t}$, $||\widetilde{g}_{t}||=e^{\beta t}$, $\langle f_{t}|e_{t} \rangle=e^{\gamma t}$ and $\langle g_{t}|e_{t} \rangle=e^{\delta t}$.

Denote the path space associated to the product system $E$ by $\Delta$. We claim that there exist functions $\phi,\psi: \Delta \to \mathbb{C}$ and a function $c:(0,\infty) \to \mathbb{C}$ which are continuous and each vanishes at the origin such that for $t>0$ and $u,v \in D(t)$,
\begin{equation}
\label{adjoint equation}
L^{e}(t,\theta_{t}u,v)=L^{e}(t,u,\theta_{t}^{*}v)+\phi(t,u)+\psi(t,v)+c(t).
\end{equation}
Here the continuity of the functions $\phi$ and $\psi$ are in the sense of Arveson. 
The desired conclusion is immediate from Eq. \ref{adjoint equation}.

Set $c(t)=(\alpha-\beta+\gamma+\delta)t$, $\phi(t,u)=L^{g}(t,u,e_t)$ and $\psi(t,v)=L^{e}(t,f_t,v)$. It is clear that $c, \phi, \psi$ are continuous and vanish at the origin. A direct calculation yields the following. For $t>0$ and $u,v \in D(t)$,
\[
e^{L^{e}(t,\theta_t u,v)}=e^{L^{e}(t,u,\theta_{t}^{*}v)+\phi(t,u)+\psi(t,v)+c(t)}.\]
Using the fact that $\phi,\psi,c$ and $L^{e}$ are continuous and vanish at the origin, we see that \[L^{e}(t,\theta_{t}u,v)=L^{e}(t,u,\theta_{t}^{*}v)+\phi(t,u)+\psi(t,v)+c(t).\] This completes the proof. \hfill $\Box$


\begin{crlre}
With the foregoing notation, we have the following.  
\begin{enumerate}
\item[(1)] The adjoint of $\widetilde{\theta}$ is $\widetilde{\theta^{*}}$.
\item[(2)] If $C_{x}$ is positive for every $x \in P$, $\widetilde{\theta}$ is positive.
\item[(3)] If $C_{x}$ is a projection for every $x \in P$, $\widetilde{\theta}$ is a projection.
\end{enumerate}
\end{crlre}
\textit{Proof.} The assertion $(1)$ is  immediate from Lemma \ref{adjoint formula}. Suppose that $C_{x}$ is positive for every $x \in P$. Note that 
$\{C_{x}^{\frac{1}{2}}\}_{x \in P}$ is a local cocycle. This is because $\theta_{x+y}^{\frac{1}{2}}=\theta_{x}^{\frac{1}{2}}\otimes \theta_{y}^{\frac{1}{2}}$ for $x,y \in P$. Let $\theta^{\frac{1}{2}}$ be the morphism associated to $\{C_{x}^{\frac{1}{2}}\}_{x \in P}$.  Note that
\[
\widetilde{\theta}=\widetilde{\theta^{\frac{1}{2}}}\widetilde{\theta^{\frac{1}{2}}}.
\]
Since $\widetilde{\theta^{\frac{1}{2}}}$ is self-adjoint, it follows that $\widetilde{\theta}$ is positive. It is routine to verify that if $\theta_{x}$  is a projection for every $x \in P$, then $\widetilde{\theta}$ is a projection. This completes the proof. \hfill $\Box$

We proceed towards computing the positive contractive local cocycles of a CCR flow. First let us fix a few notation. Let $V:=(V_{x})_{x \in P}$ be a pure isometric representation on a Hilbert space $\clh$.  For $x \in P$, we denote the range projection $V_{x}V_{x}^{*}$ by $E_{x}$ and its orthocomplement by $E_{x}^{\perp}$. The commutant of the von Neumann algebra generated by $\{V_{x},V_{x}^{*}: x \in P\}$  is denoted by $M$. By an additive cocycle of $V$, we mean a family $\xi:=\{\xi_{x}\}_{x \in P}$ of vectors in $\clh$ such that the following conditions are satisfied.
\begin{enumerate}
\item[(1)] For $x \in P$, $\xi_{x} \in Ker(V_{x}^{*})$,
\item[(2)] the map $P \ni x \to \xi_{x} \in \clh$ is continuous, and
\item[(3)] for $x,y \in P$, $\xi_{x+y}=\xi_{x}+V_{x}\xi_{y}$.
\end{enumerate}
The set of additive cocycles of $V$ is denoted by $\mathcal{A}(V)$. Note that $\mathcal{A}(V)$ is a vector space and the algebra $M$ acts on $\mathcal{A}(V)$ by $T.\xi=\{T\xi_{x}\}$ for $T \in M$ and $\xi \in \mathcal{A}(V)$. Suppose $\xi=\{\xi_{x}\}_{x \in P} \in \mathcal{A}(V)$. Note that $P \ni x \to \langle \xi_x|\xi_x \rangle \in \mathbb{R}$ is continuous and additive. Thus there exists $\mu \in \mathbb{R}^{d}$ such that $\langle \xi_x|\xi_x \rangle=\langle \mu|x\rangle$.

Let $\mathcal{D}_{0}$ be the set of triples $(\lambda,\xi,A)$ which satisfy the following conditions.
\begin{enumerate}
\item[(C1)] The element $\lambda \in \mathbb{R}^{d}$, $\xi:=\{\xi_{x}\}_{x \in P}$ is an additive cocycle of $V$ and $A$ is a positive contraction in $M$, and
\item[(C2)] for $x \in P$, $\xi_{x} \in Ker(1-A)^{\perp}$ and there exists $\eta_{x} \in Ker(1-A)^{\perp}$ such that $\xi_{x}=(1-A)^{\frac{1}{2}}\eta_x$.
\end{enumerate}
Note that  (C2) ensures that there is no ambiguity in the definition of $(1-A)^{-\frac{1}{2}}\xi_{x}$.

Let $(\lambda,\xi,A) \in \mathcal{D}_{0}$ be given. Thanks to Lemma A1 (see Appendix) of \cite{Araki}, for $x \in P$, there exists a unique  bounded linear operator $C^{(\lambda,\xi,A)}_{x}$ on $\Gamma(\clh)$ whose action on the exponential vectors is given by the formula
\[
C^{(\lambda,\xi,A)}_{x}e(\eta)=e^{-\langle \lambda|x \rangle}e^{\langle \eta|\xi_{x}\rangle}e((AE_{x}^{\perp}+E_{x})\eta+\xi_{x}). \]
Moreover Lemma A1 of \cite{Araki} implies that the norm  $||C^{(\lambda,\xi,A)}_{x}||=e^{-\langle \lambda|x \rangle}e^{||(1-A)^{-\frac{1}{2}}\xi_x||^{2}}$. Thus $C^{(\lambda,\xi,A)}_{x}$ is a contraction for every $x$ provided the following condition is satisfied.
\begin{enumerate}
\item[(C3)] For $x \in P$, $-\langle \lambda|x \rangle+||(1-A)^{-\frac{1}{2}}\xi_{x}||^{2} \leq 0$.
\end{enumerate}
Let $\mathcal{D}$ be the set of triples $(\lambda,\xi,A) \in \mathcal{D}_{0}$ for which (C3) is satisfied. 

\begin{thm}
\label{positive local cocycles}
The map \[  (\lambda,\xi,A) \to \{C_{x}^{(\lambda,\xi,A)}\}_{x \in P}\]
is a bijection between $\mathcal{D}$ and the set of positive contractive local cocycles of $\alpha^{V}$ where $\alpha^{V}$ is the CCR flow associated to $V$.
\end{thm}
\textit{Proof.}
Denote the product system associated to $\alpha^{V}$ by $E:=\{E(x)\}_{x \in P}$.
Consider an element  $(\lambda,\xi,A) \in \mathcal{D}_{0}$. Set $\widetilde{\xi}_{x}:=(1+A^{\frac{1}{2}})^{-1}\xi_{x}$. Then $\widetilde{\xi}:=\{\widetilde{\xi}_{x}\}_{x \in P}$ is an additive cocycle. Let $\mu \in \mathbb{R}^{d}$ be such that $\langle \widetilde{\xi}_{x}|\widetilde{\xi}_{x} \rangle=\langle \mu|x \rangle$. Note that $(\frac{\lambda+\mu}{2},\widetilde{\xi},A^{\frac{1}{2}}) \in \mathcal{D}_0$. A direct verification reveals that $C^{(\lambda,\xi,A)}_{x}$ is self-adjoint for every $(\lambda,\xi,A) \in \mathcal{D}_{0}$. Similarly a routine computation  leads to the equality
\[
C^{(\lambda,\xi,A)}_{x}=C^{(\frac{\lambda+\mu}{2},\widetilde{\xi},A^{\frac{1}{2}})}_{x}C^{(\frac{\lambda+\mu}{2},\widetilde{\xi},A^{\frac{1}{2}})}_{x}.\] Thus $C_{x}^{(\lambda,\xi,A)}$ is positive. We have already seen that $C_{x}^{(\lambda,\xi,A)}$ is a contraction if $(\lambda,\xi,A) \in \mathcal{D}$. 

Let $(\lambda,\xi,A) \in \mathcal{D}$ be given. 
It is clear that the map $P \ni x \to C_{x}^{(\lambda,\xi,A)}$ is weakly measurable. A routine calculation shows that $C_{x}^{(\lambda,\xi,A)} \in \alpha_{x}(B(\Gamma(\clh)))^{'}$. Let us denote $C_{x}^{(\lambda,\xi,A)}$ by $C_{x}$. For $x \in P$, let $\theta_{x}:E(x) \to E(x)$ be defined by $\theta_{x}(T)=C_{x}T$. We claim that for $x,y \in P$, $T \in E(x)$ and $S \in E(y)$,
\[
\theta_{x+y}(TS)=\theta_{x}(T)\theta_{y}(S).\]
Note that $\{T^{(x)}_{e(\xi)}: \xi \in Ker(V_{x}^{*})\}$ is total in $E(x)$. Thus, it suffices to check that for $x,y \in P$, $\xi \in Ker(V_{x}^{*})$ and $\eta \in Ker(V_{y}^{*})$, 
\[
\theta_{x+y}(T^{(x)}_{e(\xi)}T^{(y)}_{e(\eta)})=\theta_{x}(T^{(x)}_{e(\xi)})\theta_{y}(T^{(y)}_{e(\eta)}).\]
But this is straightforward. Hence $\{C_{x}^{(\lambda,\xi,A)}\}_{x \in P}$ is a positive contractive local cocycle. We leave it to the reader to make use of the purity of the representation to verify that the map 
\[  (\lambda,\xi,A) \to \{C_{x}^{(\lambda,\xi,A)}\}_{x \in P}\]
is injective. 

\textit{The surjectivity part:} Let $\{C_{x}\}_{x \in P}$ be a positive contractive local cocycle of $\alpha^{V}$. Denote the associated morphism of $E$ by $\theta:=\{\theta_x\}_{x \in P}$. Let $\widetilde{V}:=\{\widetilde{V}_{a}\}_{a \in \Omega}$ be the isometric representation constructed out of the decomposable $E_0$-semigroup $\alpha^{V}$. The Hilbert space on which $\widetilde{V}$ acts is denoted by $H_{\infty}$. In view of Prop. \ref{Isometric}, we can identify $\widetilde{V}$ with $V$ and $H_{\infty}$ with $\clh$. The identification is implemented by the unitary whose restriction to $Ker(V_{a}^{*})$, for $a \in \Omega$, is given by
\[
\xi \to [T^{(a)}_{e(\xi)}]-[T^{(a)}_{e(0)}].\]
Note that $\{\theta_{x}(T^{(x)}_{e(0)})\}_{x \in P}$ is a unit of $\alpha^{V}$. Thanks to Theorem 5.10 of \cite{Anbu}, there exist $\lambda,\mu \in \mathbb{R}^{d}$ and an additive cocycle $\xi=\{\xi_{x}\}_{x \in P}$ such that $\theta_{x}(T^{(x)}_{e(0)})=e^{-\langle \lambda|x \rangle}e^{i\langle \mu|x \rangle}T^{(x)}_{e(\xi_x)}$.

Since $\theta_{x}$ is positive, calculate as follows to observe that for $x \in P$,
\[
e^{-\langle \lambda|x \rangle}e^{i \langle \mu|x \rangle}=\langle \theta_{x}(T^{(x)}_{e(0)})|T^{(x)}_{e(0)}\rangle \geq 0.\]
This implies that $\mu=0$. 

Let $\widetilde{\theta}$ be the contraction on $H_{\infty}$ induced  by $\theta$ (See Prop. \ref{appearance of commutant}).
 Denote the contraction on $\clh$ corresponding to $\widetilde{\theta}$ by $A$. Fix $a \in \Omega$ and $\xi \in Ker(V_{a}^{*})$. Calculate as follows to observe that
\begin{align*}
[\theta_{a}T^{(a)}_{e(\xi)}]-[T^{(a)}_{e(0)}]&=[\theta_{a}T^{(a)}_{e(\xi)}]-[T^{(a)}_{e(\xi_a)}]+[T^{(a)}_{e(\xi_{a})}]-[T^{(a)}_{e(0)}] \\
                                                                        &=[\theta_{a}T^{(a)}_{e(\xi)}]-[\theta_{a}T^{(a)}_{e(0)}]+[T^{(a)}_{e(\xi_{a})}]-[T^{(a)}_{e(0)}] \\
                                                                        &=\widetilde{\theta}([T^{(a)}_{e(\xi)}]-[T^{(a)}_{e(0)}])+[T^{(a)}_{e(\xi_{a})}]-[T^{(a)}_{e(0)}] \\
                                                                        &=[T^{(a)}_{e(A\xi)}]-[T^{(a)}_{e(0)}]+[T^{(a)}_{e(\xi_{a})}]-[T^{(a)}_{e(0)}] \\
                                                                        &=[T^{(a)}_{e(A\xi+\xi_a)}]-[T^{(a)}_{e(0)}].
\end{align*}
Hence $\theta_{a}T^{(a)}_{e(\xi)}$ and $T^{(a)}_{e(A\xi+\xi_a)}$ are scalar multiplies of each other. (This is a consequence of the next Lemma whose proof we leave  to the reader.)
Thus, for $a \in \Omega$ and $\xi \in Ker(V_a^{*})$, there exists $\lambda_{a,\xi} \in \mathbb{C}\backslash \{0\}$ such that
\[
\theta_{a}T^{(a)}_{e(\xi)}=\lambda_{a,\xi}T^{(a)}_{e(A\xi+\xi_a)}.\]
Fix $a \in \Omega$ and $\xi \in Ker(V_a^{*})$. Calculate as follows to observe that
\begin{align*}
\lambda_{a,\xi}&=\langle \theta_{a}T^{(a)}_{e(\xi)}|T^{(a)}_{e(0)} \rangle \\
                        &= \langle T^{(a)}_{e(\xi)}|\theta_{a}T^{(a)}_{e(0)} \rangle \\
                        &=\langle T^{(a)}_{e(\xi)}|e^{-\langle \lambda|a \rangle}T^{(a)}_{e(\xi_a)} \rangle \\
                        &=e^{-\langle \lambda|a \rangle}e^{\langle \xi|\xi_a \rangle}.
\end{align*}
Hence we have 
\begin{equation}
\label{formula for positive cocycles}
\theta_{a}T^{(a)}_{e(\xi)}=e^{-\langle \lambda|a \rangle}e^{\langle \xi|\xi_a\rangle }T^{(a)}_{e(A\xi+\xi_a)}.
\end{equation}
Eq. \ref{formula for positive cocycles} translates to the following equation
\begin{equation}
\label{formula for contractive cocycles}
C_{a}e(\eta)=e^{-\langle \lambda|a \rangle}e^{\langle \eta|\xi_a\rangle}e((AE_a^{\perp}+E_a)\eta+\xi_a)
\end{equation}
for $a \in \Omega$ and $\eta \in \clh$. Since $C_{a}$ is bounded for every $a \in \Omega$, Lemma A1 of \cite{Araki} implies that $\xi_{a} \in Ker(1-A)^{\perp}$ and there exists $\eta_{a} \in Ker(1-A)^{\perp}$, which is necessarily unique, such that $\xi_{a}=(1-A)^{\frac{1}{2}}\eta_{a}$.

As $(1-A)^{\frac{1}{2}}$ is injective on $Ker(1-A)^{\perp}$, it is clear that $\{\eta_a\}_{a \in \Omega}$ is an additive cocycle of $V$ indexed by $\Omega$. Now apply Lemma 5.9 of \cite{Anbu} to obtain an additive cocycle indexed by $P$, which is a unique extension of $\{\eta_a\}_{a \in \Omega}$. We denote the extension by  $\{\eta_{x}\}_{x \in P}$. The density of $\Omega$ in $P$ and the fact that $\{\eta_{x}\}_{x \in P}$ and $\{\xi_{x}\}_{x \in P}$ are norm continuous imply that 
for every $x \in P$, $\xi_{x} \in Ker(1-A)^{\perp}$ and $\xi_{x}=(1-A)^{\frac{1}{2}}\eta_x$. 

Again Lemma A1 of \cite{Araki} implies that the norm of $C_{a}$, for $a \in \Omega$, is $e^{-\langle \lambda|a \rangle}e^{||\eta_{a}||^{2}}$. Since $C_{a}$ is a contraction for every $a \in \Omega$, it follows that 
\[
-\langle \lambda|a \rangle + ||\eta_{a}||^{2} \leq 0.\]
Making use of the density of $\Omega$ in $P$, we see that for every $x \in P$, 
\[
-\langle \lambda|x \rangle + ||(1-A)^{-\frac{1}{2}}\xi_{x}||^{2} \leq 0.\]
Hence $(\lambda,\xi,A) \in \mathcal{D}$. Eq. \ref{formula for contractive cocycles} implies that $C_{x}=C_{x}^{(\lambda,\xi,A)}$ for every $x \in \Omega$ and hence for every $x \in P$.  \hfill $\Box$

We have made use of the following Lemma, whose proof we leave to the reader, in the proof of Theorem \ref{positive local cocycles}.
\begin{lmma}
Let $\alpha$ be a decomposable $E_0$-semigroup with product system $E:=\{E(x)\}_{x \in P}$. Fix a left coherent section $(e_x)_{x \in P}$ of decomposable vectors with unit norm. Denote the isometric representation built out of $\alpha$ by $V$ and let $H_{\infty}$ be the Hilbert space on which $V$ acts. Let $a \in \Omega$ and $u,v \in D(a)$ be given. Suppose 
\[
[u]-[e_a]=[v]-[e_a]
\]
then $u$ and $v$ are scalar multiples of each other.

\end{lmma}

It is relatively easy to write down the set of projective local cocycles of $\alpha^{V}$ from Theorem \ref{positive local cocycles}. We merely state the result and omit the details. Let $\mathcal{E}$ be the set of pairs $(\xi,Q)$ where 
\begin{enumerate}
\item[(E1)] the operator $Q$ is a projection in $M$ and $\xi:=\{\xi_{x}\}_{x \in P}$ is an additive cocycle of $V$, and
\item[(E2)] for every $x \in P$, $(1-Q)\xi_{x}=\xi_{x}$.
\end{enumerate}
Consider an element $(\xi,Q) \in \mathcal{E}$. Then there exists a unique $\lambda_{\xi} \in \mathbb{R}^{d}$ such that for $x \in P$, $\langle \lambda_{\xi}|x \rangle=||\xi_x||^{2}$. 

\begin{ppsn}
The map 
\[
(\xi,Q) \to \{C_{x}^{(\lambda_{\xi}, \xi,Q)}\}_{x \in P}\]
is a bijection between $\mathcal{E}$ and the set of projective local cocycles of $\alpha^{V}$.
\end{ppsn} 

\begin{rmrk}
\label{commutant}
Suppose the dimension of $P$ is at least $2$. 
For the isometric representations that arise out of $P$-modules (See Eq. \ref{isometries}) the additive cocycles and the commutant of the representation were computed in \cite{Anbu_Sundar}. 
Let $A$ be a $P$-module and $V^{(A)}$ be the isometric representation associated to $A$ on $L^{2}(A)$. It was shown in \cite{Anbu_Sundar} that $V^{(A)}$ admits no nontrivial additive cocycle. 
(See Prop. 2.4 of \cite{Anbu_Sundar}).
Let \[
G_{A}=\{z \in \mathbb{R}^{d}: A+z=A\}.
\]
Then $G_{A}$ is a closed subgroup of $\mathbb{R}^{d}$ called the isotropy group of $A$. Note that $G_{A}$ acts on $L^{2}(A)$ by translations. One of the main result in \cite{Anbu_Sundar} is that
the commutant of $V^{(A)}$ is generated by $G_{A}$. (See Corollary 3.4 of \cite{Anbu_Sundar}).


\end{rmrk}

\section{Prime CCR flows}
As another application, we derive a necessary and a sufficient condition for a  CCR flow to be prime. We start with the following definition.

\begin{dfn}
Let $\alpha:=\{\alpha_x\}_{x \in P}$ be an $E_0$-semigroup. We say that $\alpha$ is prime if whenever $\alpha$ is cocycle conjugate to $\beta \otimes \gamma$ where $\beta$ and $\gamma$ are $E_0$-semigroups then
either $\beta$ or $\gamma$ is an automorphism group. 
\end{dfn}

The problem of constructing prime $E_0$-semigroups, in the $1$-parameter case, has recieved considerable attention in the recent years. We refer the reader to the paper \cite{Jankowski} and the references therein for more details. In the one 
parameter case, the interesting question is to produce prime $E_0$-semigroups which are not of type I. The reason being that we do know which type I $E_0$-semigroups or equivalently which CCR flows are prime. From the classification of type I $E_0$-semigroups and on an index computation, it is well known
that a $1$-parameter CCR flow is prime if and only if the ``pure" part of the corresponding isometric representation is unitarily equivalent to the standard shift semigroup on $L^{2}(0,\infty)$. 

It is only natural to ask which CCR flows in the multiparameter case are prime. The main aim of this section is to prove the following theorem. We need a  bit of terminology. Let $V:=(V_{x})_{x \in P}$ be an isometric representation on a Hilbert space $\clh$. We say that $V$ is irreducible if the only closed subspaces of $\clh$ invariant under $\{V_{x}, V_{x}^{*}:x \in P\}$ are $\{0\}$ and $\clh$.

\begin{thm}
\label{prime}
Let $V:=\{V_{x}\}_{x \in P}$ be a pure isometric representation. Denote the CCR flow corresponding to $V$ by $\alpha^{V}$.
The following are equivalent.
\begin{enumerate}
\item[(1)] The CCR flow $\alpha^{V}$ is prime.
\item[(2)] The representation $V$ is irreducible.
\end{enumerate}
\end{thm}
Note that the CCR functor converts direct sum of isometric representations into tensor product. Thus the implication $(1) \implies (2)$ is obvious. 
The first lemma that we need is the following. 

Let $\alpha^{(1)},\alpha^{(2)}$ be $E_0$-semigroups and set $\alpha:=\alpha^{(1)} \otimes \alpha^{(2)}$.  Let $E^{(i)}:=\{E^{(i)}(x)\}_{x \in P}$ be the product system associated to
$\alpha^{(i)}$. Similarly let $E:=\{E(x)\}_{x \in P}$ be the product system associated to $\alpha$. For $x \in P$, we denote the set of decomposable vectors in $E(x)$, $E^{(i)}(x)$ by $D(x)$ 
and $D^{(i)}(x)$ respectively. With the foregoing notation, we have the following.

\begin{lmma}
Suppose $\alpha$ is decomposable. Then $\alpha^{(1)}$ and $\alpha^{(2)}$ are decomposable. Moreover for $x \in P$,
we have
\[
D(x)=\{u_1 \otimes u_2: u_1 \in D^{(1)}(x), u_2 \in D^{(2)}(x)\}.\]
\end{lmma}
\textit{Proof.} It is clear that for $x \in P$, $\{u_1 \otimes u_2: u_1 \in D^{(1)}(x), u_2 \in D^{(2)}(x)\} \subset D(x)$.  Suppose $P=[0,\infty)$. From Arveson's characterisation of decomposable product systems and type I product systems, it follows that $\alpha^{(1)}$ and $\alpha^{(2)}$ are CCR flows.
The equality of decomposable vectors follows from Prop. \ref{decomposable}. Thus the conclusion is true in the one parameter case. 

Now let $d \geq 2$. Fix $x \in P$. By restricting to the ray $\{tx: t \geq 0\}$ and from the one parameter conclusion, we see that if $u \in D(x)$ then there exists $u_1 \in E^{(1)}(x)$ and $u_2 \in E^{(2)}(x)$ such that $u=u_1 \otimes u_2$. Thus, for $x \in P$, every decomposable vector of $E(x)$ splits us a product $u_1 \otimes u_2$ with $u_i \in E^{(i)}(x)$. We leave it to  the reader to make use of Remark \ref{A few remarks} to convince herself/himself that this precisely implies that
\[
D(x)=\{u_1 \otimes u_2: u_1 \in D^{(1)}(x), u_2 \in D^{(2)}(x)\}.\]
The rest of the conclusion is immediate and the proof is complete. \hfill $\Box$

Keep the foregoing notation. Suppose that $\alpha$ is decomposable. For $i=1,2$, let $e^{(i)}:=(e^{(i)}_{x})_{x \in P}$ be a left coherent section of decomposable vectors of $E^{(i)}$ of unit norm. For $x \in P$, set $e_x:=e^{(1)}_{x} \otimes e^{(2)}_{x}$. Note that $e:=(e_x)_{x \in P}$ is a left coherent section of $E$ consisting of decomposable vectors with unit norm. It follows immediately from the definition that 
for $a \in \Omega$, $u^{(i)} \in D^{(i)}(a)$ and $v^{(i)} \in D^{(i)}(a)$ 
\begin{equation}
\label{product of path spaces}
L^{e}(a,u^{(1)}\otimes u^{(2)},v^{(1)}\otimes v^{(2)})=L^{e^{(1)}}(a,u^{(1)},v^{(1)})+L^{e^{(2)}}(a,u^{(2)},v^{(2)}).
\end{equation}
An immediate consequence of Eq. \ref{product of path spaces} is that the isometric representation constructed out of $\alpha$ is the direct sum of isometric representations constructed out of $\alpha^{(i)}$. Let us explain this briefly.
Let $V^{(i)}$ be the isometric representation constructed out of $\alpha^{(i)}$ which acts on $H^{(i)}_{\infty}$ and $V$ be the isometric representation constructed out of $\alpha$ which acts on $H_{\infty}$. 
Then the map
\[
H^{(1)}_{\infty} \oplus H^{(2)}_{\infty} \ni ([u^{(1)}]-[v^{(1)}])\oplus([u^{(2)}]-[v^{(2)}]) \to [u^{(1)}\otimes u^{(2)}]-[v^{(1)}\otimes v^{(2)}] \in H_{\infty}
\]
is a unitary which intertwines $V^{(1)}\oplus V^{(2)}$ and $V$.

We need one more little lemma before we can prove Theorem \ref{prime}.

\begin{lmma}
\label{automorphism}
Let $\alpha$ be a decomposable $E_0$-semigroup with product system $E$. The isometric representation constructed out of $\alpha$ will be denoted by $V$ and the Hilbert space on which it acts will be denoted $H_{\infty}$. Suppose that $H_{\infty}=\{0\}$. Then $\alpha$ is an automorphism group, i.e. for every $x \in P$, $\alpha_{x}$ is an automorphism.
\end{lmma}
\textit{Proof.} Let $e:=(e_x)_{x \in P}$ be a left coherent section of decomposable vectors with unit norm. Consider an element $a \in \Omega$. Suppose $u \in D(a)$ be such that $\langle u|e_a \rangle=1$. 
The equality $\langle [u]-[e_a]|[u]-[e_a] \rangle=0$ implies that $L^{e}(a,u,u)=0$. Taking exponential, we get $\langle u|u \rangle=1$. Thus we get the equality
\[|\langle u|e_{a}\rangle| = ||u||||e_{a}||.\]
Hence $u$ is a scalar multiple of $e_{a}$. But $D(a)$ is total in $E(a)$. This has the consequence that $E(a)$ is one-dimensional for every $a \in \Omega$. 

For $x \in P$, let $d(x)=\dim E(x)$. We have the equality $d(x+a)=d(x)d(a)$ for every $x \in P$ and $a \in \Omega$. But $\Omega$ is an ideal in $P$ and $d(a)=1$ for every $a \in \Omega$. Hence $d(x)=1$ for every $x \in P$.
In other words, $E(x)$ is one-dimensional or equivalently $\alpha_{x}$ is an automorphism for every $x \in P$. This completes the proof. \hfill $\Box$

\textit{Proof of Theorem \ref{prime}.} We have already proved the implication $(1) \implies (2)$. Assume that $(2)$ holds. Suppose that $\alpha^{V}=\alpha^{(1)}\otimes \alpha^{(2)}$. For $i=1,2$, let $V^{(i)}$ be the isometric representation corresponding to $\alpha^{(i)}$. By our preceeding discussions, it follows that $V$ is equivalent to the direct sum $V^{(1)}\oplus V^{(2)}$. But $V$ is irreducible. Hence either $V^{(1)}$ or $V^{(2)}$ must be zero. Lemma \ref{automorphism} implies that either $\alpha^{(1)}$ or $\alpha^{(2)}$ is an automorphism group. This shows that $\alpha^{V}$ is prime. The proof is now complete. \hfill $\Box$

\begin{rmrk}
It is a difficult problem to determine all irreducible isometric representations of an arbitrary cone. Up to the author's knowledge, this question remains open even for the simplest case of the quarter plane. However what we do know is that when the dimension of the cone, i.e. $d \geq 2$, there are indeed uncountably many irreducible isometric representations. Consequently, there are uncountably many prime CCR flows in the multiparameter case.

Suppose that $d \geq 2$. Taking into account Remark \ref{commutant} and the discussions following Theorem \ref{Injectivity of CCR}, it suffices to exhibit an uncountable family $\mathcal{F}$ of $P$-modules 
such that
\begin{enumerate}
\item[(1)] if $A \in \mathcal{F}$, the isotropy group $G_{A}$ is trivial, and
\item[(2)] if $A,B \in \mathcal{F}$ with $A \neq B$ then $A$ is not a translate of $B$.
\end{enumerate}
Choose $a \notin P \cup -P$. 
 The family  $\{P \cup (ta+P): t>0\}$ is one such candidate.
\end{rmrk}

\bibliography{references}
 \bibliographystyle{amsplain}
 
 \nocite{Araki}
 
\noindent
{\sc S. Sundar}
(\texttt{sundarsobers@gmail.com})\\
         {\footnotesize  Institute of Mathematical Sciences, CIT Campus, \\
Taramani, Chennai, 600113, Tamilnadu, INDIA.}\\

\end{document}